\documentclass[12pt]{amsart}
\usepackage{amsthm}
\usepackage{amscd}
\usepackage{calrsfs}
\usepackage{url}
\usepackage{graphicx}
\usepackage{tikz-cd}
\usepackage{mathtools}
\usepackage{color}
\usepackage{amsmath}
\usepackage{amssymb}
\usepackage{amsthm}
\usepackage{amscd}
\usepackage{eucal} 
\usepackage{enumitem}
\usepackage{mathdots}
\usepackage{tikz}
\usepackage{hyperref}
\usepackage[all]{xy}

\allowdisplaybreaks

\newtheorem{theorem}[equation]{Theorem}

\newtheorem{lemma}[equation]{Lemma}

\newtheorem{proposition}[equation]{Proposition}
\newtheorem{corollary}[equation]{Corollary}

\newtheorem*{theorem:derhamisomorphism}{Theorem~\ref{T:derhamisomorphism}}
\newtheorem*{theorem:characterization}{Theorem~\ref{T:characterization}}

\theoremstyle{definition}
\newtheorem{definition}[equation]{Definition}

\newtheorem{remark}[equation]{Remark}

\numberwithin{equation}{section}

\newcommand{\ZZ}{\mathbb{Z}}

\newcommand{\TT}{\mathbb{T}}

\newcommand{\CC}{\mathbb{C}}
\newcommand{\C}{\mathbb{C}}

\newcommand{\bA}{\mathbf{A}}
\newcommand{\bV}{\mathbf{V}}

\newcommand{\bm}{\mathbf{m}}
\newcommand{\bn}{\mathbf{n}}

\newcommand{\bv}{\mathbf{v}}
\newcommand{\bz}{\mathbf{z}}

\newcommand{\bb}{\mathbf{b}}

\newcommand{\cM}{\mathcal{M}}

\newcommand{\zz}{\mathbf{z}}

\DeclareMathOperator{\Log}{Log}

\DeclareMathOperator{\GL}{GL}

\DeclareMathOperator{\Lie}{Lie}
\DeclareMathOperator{\Mat}{Mat}

\DeclareMathOperator{\ch}{ch}

\DeclareMathOperator{\Frob}{Frob}

\DeclareMathOperator{\Exp}{Exp}
\DeclareMathOperator{\diag}{diag}

\DeclareMathOperator{\Gal}{Gal}

\DeclareMathOperator{\dmat}{d}

\DeclareMathOperator{\Res}{Res}

\DeclareMathOperator{\Id}{Id}

\newcommand{\tr}{\mathrm{tr}}
\newcommand{\tpi}{\widetilde{\pi}}

\definecolor{ForestGreen}{rgb}{0.0, 0.5, 0.0}

\begin{document}

\title[Mellin Transform formulas ]{Mellin Transform formulas for special values of Goss $L$-functions}

	\author{O\u{g}uz Gezm\.{i}\c{s}}
\address{Department of Mathematics, National Tsing Hua University, Hsinchu City 30042, Taiwan R.O.C.}
\email{gezmis@math.nthu.edu.tw}


\author{Nathan Green}

\address{Department of Mathematics and Statistics, University of North Carolina at Charlotte, Charlotte, NC 28223,  U.S.A.}
\email{ngreen50@charlotte.edu}

\date{\today}



    \begin{abstract} In the present paper, we provide formulas for the special values of dual Goss $L$-functions of Drinfeld modules at any non-negative integer in terms of the determinant of the values of a motivic map extensively studied in \cite{GG26} for a large class of Anderson modules. Our result may be seen as a function field analogue of conjectural formulas for special values of motivic $L$-functions in terms of algebraic Mellin integrals.
    \end{abstract}
    \maketitle
    
\section{Introduction}
\subsection{Motivation and Background}
Ever since Euler discovered the surprising result that
\[\zeta(2) = \frac{\pi^2}{6},\]
mathematicians have sought to extend and understand this connection between periods and zeta values. On the one hand, mathematicians developed an entire framework related to periods that extends $\frac{\pi^2}{6}$ to other meaningful arithmetic quantities. The work of Deligne, Grothendieck, Kontsevich, Zagier and others (see \cite{KZ01} and the references therein) demonstrated that $\pi$ and other related transcendental quantities can be viewed as specific examples of a large class of numbers called \textit{periods} (see Remark \ref{R:KZperiods} for their algebro-geometric definition). These are values that result from a natural pairing between the Betti and de Rham cohomology of certain varieties. One can view this as a natural generalization of the fact that
\begin{equation}\label{E:periodclassical}
2\pi i = \oint_\sigma \frac{dz}{z},
\end{equation}
where $\sigma$ is a counter-clockwise unit circle path about the origin.

On the other hand, numerous mathematicians pioneered various settings for generalizing the Riemann zeta function. Initially, this included Dedekind zeta functions, Dirichlet $L$-functions and other straightforward generalizations of the Riemann zeta function. Each of these constructions was connected back to various periods or other meaningful arithmetic data, such as in the class number formula. Eventually, this led to the framework of $L$-functions for various geometrical objects, such as $L$-functions of elliptic curves and abelian varieties. This led Deligne, Grothendieck, Serre and others (see \cite{Del} and the references therein) to construct the $L$-function of a motive. Motives are algebraic-geometric objects that live in a conjectural category which generalizes the various cohomology theories attached to a variety.

Building on these two paths of generalizations, mathematicians, including Beilinson \cite{Bei}, Bloch and Kato \cite{BK}, and Deligne \cite{Del}, formulated conjectures connecting the two sides of this world; relating special values of motivic $L$-functions with periods of motives. Specifically, for a pure motive $M$ defined over $\mathbb{Q}$, in \cite{Del}, Deligne defined quantities $c^\pm(M)$ which are the determinants of certain comparison isomorphisms between variants of Betti and de Rham cohomology. It is possible to choose bases for these cohomology spaces, then to write these quantities $c^\pm(M)$ explicitly as the determinants of a matrix of  periods (see in particular \cite[(1.7.2)]{Del} and the following discussion therein).

More recently, and more relevant to the results of this paper, in \cite{Bro16}, Brown described how to give a general geometric framework for the integrals that give rise to the periods described above. They are given as period integrals of moduli spaces of curves of genus zero with $n$ ordered marked points. In 2026, Brown \cite{Bro26} observed that such integrals are a special case of the more general class of integrals called \textit{algebraic Mellin integrals}. Specifically, an algebraic Mellin integral on a smooth affine variety $X$ of dimension $d$ is of the form
\begin{equation}\label{D:alg Mellin}
I(s_1,\dots,s_r) = \int_\sigma f_1^{s_1} \dots f_r^{s_r} \omega,
\end{equation}
where $s_i\in \C$ with nonnegative real part, $\omega$ is a regular differential $d$ form on $X$,  $\sigma \subset X(\C)$ is a relative $d$-chain with a certain boundary condition and each $f_i$ is a regular function on $X$. This can be interpreted as a generalization of the pairing of Betti and de Rham cohomology discussed above, with $\sigma$ coming from some relative Betti homology class and $\omega$ from the de Rham cohomology. Brown then goes on to construct a matrix of such algebraic Mellin integrals and studies its determinant for transcendence considerations.

In this paper, we study similar problems in the setting of characteristic $p$ function fields. We construct matrices of periods that arise from a Drinfeld module analogue of Mellin integrals and show that the determinant of this matrix is equal to a special value of the dual Goss $L$-function. There have been some results in this direction given previously. For example, Fang \cite{Fang} showed that the determinant of a function field `regulator' (a matrix of logarithm values) gives Taelman $L$-values. Additionally, Angl\'es, Ngo Dac and Tavares Ribeiro showed \cite{ANDTR20} (using different techniques) that Taelman $L$-values are, up to a rational constant, given by the determinant of a certain space generated by logarithms. The novelty of our approach is that it relates dual Goss $L$-values directly to a determinant of Drinfeld module Mellin integrals. The bulk of the delicate analysis for our formulas in this paper was carried out in \cite{GG26}, where the authors proved formulas for special values of dual Goss $L$-functions only at $s=0$. In this paper, we extend those techniques and carry out the analysis so that we can give formulas of $L$-functions evaluated at arbitrary positive integers. Further, our formulas motivate us to study a function field analogue of Kontsevich-Zagier periods and we briefly discuss, in Remark \ref{R:KZperiods}, a future research direction aiming to explicitly construct these objects in our setting.

\subsection{Summary of Main Results} \label{S:mainresults}We now provide a sufficient background in order to state our main results and explain their context. Let $\mathbb{F}_q$ be a finite field with $q$ elements. Letting $\theta$ be a variable over $\mathbb{F}_q$, we define $A$ to be the ring of polynomials in $\theta$ with coefficients in $\mathbb{F}_q$. We set $K$ to be the fraction field of $A$. We consider the norm $|\cdot|$ corresponding to the infinite place normalized so that $|\theta|=q$. We let $K_{\infty}$ be the completion of $K$ with respect to $|\cdot|$ which is indeed given by the formal Laurent series ring $\mathbb{F}_q((1/\theta))$. We also let $\mathbb{C}_{\infty}$ the completion of a fixed algebraic closure of $K_{\infty}$ and $\overline{K}$ be a fixed algebraic closure of $K$ in $\mathbb{C}_{\infty}$.

Let  $t$ be a variable over $\mathbb{C}_{\infty}$. We define \textit{the Tate algebra $\mathbb{T}$} by 
\[
\mathbb{T}:=\left\{\sum_{i=0}^{\infty}a_it^i\in \mathbb{C}_{\infty}[[t]] \ \ | \ \ |a_i|\to 0 \text{ as } i\to \infty \right\}.
\]
It is equipped with the \textit{Gauss norm} $||\cdot ||$ defined by 
\[
||g||:=\max\{|a_i| \ \ | \ \ i\in \mathbb{Z}_{\geq 0}\}, \ \ \ \ \ \ \ g=\sum_{i\geq 0}a_it^i\in \mathbb{T}.
\]
For any $j\in \mathbb{Z}$, we also define \textit{the $j$-th twisting operator} $(j):\mathbb{T}\to \mathbb{T}$ by
\[
g^{(j)}:=\sum_{i=0}^{\infty}a_i^{q^j}t^i\in \mathbb{T}, \ \ \ \ \ \ g=\sum_{i=0}^{\infty}a_it^i\in \mathbb{T}.
\]
For any $G=(g_{\mu \nu})_{\mu \nu}\in \Mat_{d\times \ell}(\mathbb{T})$, by a slight abuse of notation, we further define
$
||G||:=\max\{||g_{\mu \nu}|| \ \ | \ \ 1\leq \mu, \nu \leq r\}
$. When $G\in \Mat_{d\times \ell}(\mathbb{C}_{\infty})$, we use the earlier notation and let $|G|:=||G||$.  We further set 
\[
G^{(j)}:=(g_{\mu \nu}^{(j)})_{\mu \nu}\in \Mat_{d\times \ell}(\mathbb{T}).
\]
Let $L$ be an $A$-algebra. We denote by $ \Mat_{d}(L[[\tau]]) $ the non-commutative power series ring and by $ \Mat_{d}(L[\tau]) $ the non-commutative polynomial ring  which is defined subject to the condition
\[
\tau B=B^{(1)} \tau, \ \ B\in \Mat_{d}(L).
\]
There exists an action of $ \Mat_{d}(L[[\tau]]) $ on $L^d:=\Mat_{d\times 1}(L)$ defined by
\[
\mathcal{U}\cdot x=\sum_{i=0}^{\infty}B_ix^{(i)}, \ \ x\in L^d, \ \ \ \ \mathcal{U}=\sum_{i\geq 0}B_i\tau^i\in \Mat_{d}(L[[\tau]])
\]
whenever the right hand side of the first equality above converges in $L$.

Let $\bA:=\mathbb{F}_q[t]$. By an abelian Anderson module of dimension $d$ defined over a field extension $L$ of $K$, we mean an $\mathbb{F}_q$-algebra homomorphism $\phi:\bA\to \Mat_d(L[\tau])$ given by 
       \[
       \phi_{t}=\dmat_{\phi}[t]+A_1\tau+\cdots+A_{s}\tau^s
       \]
       for some $s\in \ZZ_{\geq 0}$ so that $\dmat_{\phi}[t]:=\theta \Id_d+N$ where $N\in \Mat_d(L)$ is a nilpotent matrix. For each $\phi$, there exists a uniquely defined everywhere convergent $\mathbb{F}_q$-linear map $\Exp_{\phi}:\mathbb{C}_{\infty}^d\to \mathbb{C}_{\infty}^d$, which is called \textit{the exponential function of $\phi$}. Its kernel is a free and finitely generated $\bA$-module and is called \textit{the period lattice of $\phi$}. We call its non-zero elements \textit{periods of $\phi$} (see \S\ref{S:Abelian And Modules} for more details). Furthermore, we call $\phi$ \textit{uniformizable} if $\Exp_{\phi}$ is a surjective map. One can define \textit{the logarithm function} $\Log_{\phi}:\mathcal{D}\to \mathbb{C}_{\infty}^{d}$ which is the inverse of $\Exp_{\phi}$ converging on certain subdomain $\mathcal{D}$ of $\mathbb{C}_{\infty}^d$. Moreover, we note that one can also naturally extend the domain of $\Exp_{\phi}$ to $\mathbb{T}^d=\Mat_{d\times 1}(\mathbb{T})$ and the domain of $\Log_{\phi}$ to a certain subset of $\mathbb{T}^d$ (see \S\ref{S:Abelian And Modules} for more details on this extension).
       
Let $\phi$ be an abelian and uniformizable Anderson module. Let $\{\lambda_1,\dots,\lambda_r\}$ be a generating set of the period lattice of $\phi$ as an $\bA$-module and let $\zz\in \mathbb{C}_{\infty}^d$. For $1\leq \nu \leq r$, consider the Anderson generating function of $\phi$ with respect to $\lambda_{\nu}$ (see \S\ref{S:Abelian And Modules} for more details on Anderson generating functions) given by
\[
\mathcal{Y}_{\nu}(t):=\sum_{i=0}^{\infty}\Exp_{\phi}\left(\dmat_{\phi}[t]^{-(i+1)}\lambda_{\nu}\right)t^i=(\mathcal{Y}_{\nu 1},\dots,\mathcal{Y}_{\nu d})^{\tr}\in \mathbb{T}^d.
\]
    We define 
\begin{multline}\label{D:M phi}
\mathbb{M}_{\phi}:=\{\left(\sum_{n=0}^{\infty}a_{1,n}\tau^n,\dots,\sum_{n=0}^{\infty}a_{d,n}\tau^n\right)\in \Mat_{1\times d}(\mathbb{C}_{\infty}[[\tau]])\ \ | \ \  \\ \text{ for each } 1\leq \ell \leq d\ \ \text{ and }1\leq \nu \leq r, \ \  
|a_{\ell,n}|||\mathcal{Y}_{\nu \ell}^{(n)}||\to 0 
\text{ as } n\to \infty \}.
\end{multline}
We equip $\mathbb{M}_{\phi}$ with the norm $|\cdot|_{\phi}$ defined by 
\[
|\mathcal{G}|_{\phi}:=\max\{|a_{\ell,n} | ||\mathcal{Y}_{\nu\ell}^{(n)}|| \ \ | \ \ n\geq 0,\ \ \ 1\leq \ell \leq d,  \ \ 1\leq \nu \leq r\}.
\]
Observe that $(\mathbb{M}_{\phi},|\cdot|_{\phi})$ is a normed $\mathbb{C}_{\infty}$-vector space. 
In \cite[\S3.2]{GG26},  we define the map
$\varphi_{\phi}:\mathbb{M}_{\phi}\to \Mat_{1\times r}(\mathbb{T})$ given by   
\begin{equation}\label{E:var phi intro}
\varphi_{\phi}\left(\left(\sum_{n=0}^{\infty}a_{1,n}\tau^n,\dots,\sum_{n=0}^{\infty}a_{d,n}\tau^n\right)\right):=\sum_{1\leq \ell \leq d}\sum_{n=0}^{\infty}a_{\ell,n}\sum_{\mu=1}^{r}(\mathcal{B}_{\mu}^{[\ell]})^{(n)}\mathfrak{e}_{\mu}^{\tr}\mathfrak{S}_{n-1}.
\end{equation}
Here, $  \mathcal{B}_{\mu}^{[\ell]}\in \mathbb{C}_{\infty}[t]$, $\mathfrak{e}_{\mu}\in \Mat_{r\times 1}(\mathbb{F}_q)$ is the $\mu$-th unit vector and $\mathfrak{S}_{n-1}\in \Mat_r(\mathbb{C}_{\infty}[t])$ (see \S \ref{S:map phi} for precise definitions). It is shown in  \cite[Thm. 3.7]{GG26} that $\varphi_\phi$ is well-defined, continuous and injective.

Let $M=\Mat_{1\times d}(\mathbb{C}_{\infty}[\tau])$ be the $t$-motive of $\phi$ (see \S\ref{S:tmotive} for more details on $t$-motives). Note that, via our map $\varphi_{\phi}$, we have $M\cong \Mat_{1\times r}(\mathbb{C}_{\infty}[t])$ as $\mathbb{C}_{\infty}[t,\tau]$-modules. For a fixed $\bz\in \C_\infty^d$, we define $\delta_{1,\zz}^M:M\to \mathbb{C}_{\infty}$ by 
\[
\delta_{1,\zz}^M(m):=m\zz:=m_1\cdot z_1+\dots+m_d \cdot z_d, \ \ m=[m_1,\dots,m_d]\in \Mat_{1\times d}(\mathbb{C}_{\infty}[\tau]).
\]
 We also define $M_{\zz}$ to be the set of elements $(\mathfrak{c}_1,\dots,\mathfrak{c}_d)\in \Mat_{1\times d}(\mathbb{C}_{\infty}[[\tau]]) $ where, for each $i\in \{1,\dots,d\}$ and $\mathfrak{c}_i=\sum_{j=0}^{\infty}c_{i,j}\tau^j\in \mathbb{C}_{\infty}[[\tau]]$, we have
$
\left(c_{1,\mu}\tau^\mu,\dots,c_{d,\mu}\tau^\mu \right)\zz\to 0 $ as $\mu\to \infty$. Then we extend the map $\delta_{1,\zz}^M$ to $M_{\zz}$ by defining $\delta_{1,\zz}^{M}:M_{\zz}\to \mathbb{C}_{\infty}$ to be
\[
\delta_{1,\zz}^M(\tilde{m}):=\lim_{\mu\to \infty}\delta_{1,\zz}([\mathfrak{c}^{\mu}_1,\dots,\mathfrak{c}^{\mu}_d])
\]
where
$
\tilde{m}=\left[\mathfrak{c}_1,\dots, \mathfrak{c}_d\right]
$
and 
$\mathfrak{c}_i^{\mu}:=\sum_{j=0}^{\mu}c_{i,j}\tau^j$. 
 Finally, we define the map $\mathcal{M}_{\zz}:K\subset \Mat_{1\times r}(\mathbb{T}) \to \C_\infty$ given by 
\[
\mathcal{M}_{\zz}:=\delta_{1,\zz}^{M}\circ \varphi_{\phi}^{-1}.
\]
Technically the domain of $\mathcal{M}_{\zz}$ is $\varphi_{\phi}(M_\zz\cap\mathbb M_\phi)$, but we do not emphasize this here, since we will check well-definedness manually each time we evaluate an element under this map).  

We view the map $\mathcal{M}_{\zz}$ as a Drinfeld module analogue of the Mellin transform. The construction of $\mathcal{M}_{\zz}$ goes back to \cite{Gre24} of the second author. There, in the case of tensor powers of the Carlitz module (which are akin to Tate twists of the multiplicative group) the second author defined a motivic operator $\cM_\bz$ that mimics certain aspects of the classical Mellin transform. He showed that using this operator, we get (as a special case when $n=1$)
\[
\cM_{\bz}\left(\frac{u}{\exp_C(u)}\right ) = \Gamma_A(n)\zeta_A(n) \in K_\infty,
\]
where $u = \frac{\tpi}{\theta-t}$ ($\tpi$ is the Carlitz period \eqref{E:tildepi}) and $\Gamma_A$ and $\zeta_A$ are function field analogues of the gamma and Riemann zeta functions, respectively. The above equation is a special case of \cite[Cor. 5.11]{Gre24} for $n=1$; it is proved for general $n$ in \cite{Gre24}. The reader should compare this with the formula using the classic Mellin transform,
\[\cM\left(\frac{1}{e^x-1}\right )(s) = \Gamma(s)\zeta(s),
\]
where $\cM$ is the Mellin transform (see \cite[(1.1)]{Gre24}) and $\Gamma(s)\zeta(s)$ is the completed zeta function.

Let $\mathfrak{N}=(m_1,\dots,m_d)\in \mathbb{M}_{\phi}$. For any $g=(g_1,\dots,g_d)^{\tr}\in \mathbb{T}^d$ such that $\lVert g_i \rVert \leq \lVert \mathcal{Y}_{\nu i} \rVert$ for each $1\leq \nu \leq r$, we define 
 \[
 \mathfrak{N}(g):=m_1\cdot g_1+\cdots+m_d\cdot g_d\in \TT.
 \]
We also set
\begin{equation}\label{D:mathcal F}
\mathcal{F}:=[\mathcal{Y}_{1}(t),\dots,\mathcal{Y}_{r}(t)]\in \Mat_{d\times r}(\mathbb{T}).
\end{equation}
We let $\Psi\in \GL_r(\TT)$ be the matrix defined as in \cite[\S2.4.1]{GG26}, which is \textit{a rigid analytic trivialization of the $t$-comotive of $\phi$} (see \S\ref{S:Drinfeld}, \S\ref{S:Carlitztensorcase} and \S\ref{S:drinfeldtens} for the construction of $\Psi$ when $\phi$ is a particular Anderson module).  We finally set (for $\zz$ sufficiently small)
\[ \mathcal{M}_{\zz}(\mathfrak{N}\cdot \mathcal{F}((\Psi)^{\tr})^{(-1)}):=\mathcal{M}_{\zz}([\mathfrak{N}(\mathcal{Y}_1),\dots,\mathfrak{N}(\mathcal{Y}_r)]((\Psi)^{\tr})^{(-1)}).
 \]
We comment that by \cite[Prop. 3.6]{GG26}, each $[\mathfrak{N}(\mathcal{Y}_1),\dots,\mathfrak{N}(\mathcal{Y}_r)]((\Psi)^{\tr})^{(-1)}$ is in the range of $\varphi_{\phi}$, namely it is equal to $\varphi_{\phi}(\mathfrak{N})$. Thus there is no issue in taking $\varphi^{-1}_{\phi}$ in the above formula. Given an  Anderson module $\phi$, we also remind the reader of the construction of the non-negative integer $\mathfrak{m}_\phi$ described in \cite[Lem. 3.5]{GG26}, that is the smallest integer which causes $1/\theta^{\mathfrak{m}_\phi} \Exp_{\phi}\left(\dmat_{\phi}[t]^{-(i+1)}\lambda_{\nu}\right)$ to be inside the domain of convergence of the logarithm function of $\phi$.

\textit{A Drinfeld module $\phi$ of rank $r\geq 1$ defined over a field extension $L$ of $K$} is a one dimensional Anderson module, namely it is an $\mathbb{F}_q$-algebra homomorphism $\phi:\bA\to L[\tau]$ given by
\begin{equation}\label{E:Drinfeld}
\phi_{t}=\theta+k_1\tau+\cdots+k_r\tau^r\in L[\tau], \ \ k_r\neq 0.
\end{equation}
Let $\rho^{[k]}$ be the abelian Anderson module corresponding to the tensor product $M_{\phi}\otimes_{L[t]} M_{C^{\otimes k}}$ of $t$-motives introduced in \S\ref{S:drinfeldtens}. By convention, when $k=0$, we set $\rho^{[k]}=\phi$.  For $\phi = \rho^{[k]}$, we will shorten our notation by setting $\mathfrak{m}_k:=\mathfrak{m}_{\rho^{[k]}} $.

For each $1\leq \ell \leq rk+1$ and $k\geq 1$, we set $\mathfrak{N}_{\ell}\in \mathbb{M}_{\rho^{[k]}}$ to be the $\ell$-th row of $\Log_{\rho^{[k]}}\frac{1}{\theta^{\mathfrak{m}_k}}$. Let $\omega_C\in \mathbb{T}$ be the Anderson-Thakur function defined in \eqref{D:omega_C} and $f_{\lambda_i}$ and $\Upsilon_{\phi}\in \GL_r(\mathbb{T})$ be as given in \eqref{E:Upsilon}. For each $k\geq 0$, we define the block matrix
\[
\mathcal{F}_k:=\omega_C^k\begin{pmatrix}
    & \Upsilon_{\phi} &\\
    & (t-\theta)\Upsilon_{\phi} & \\
    & \vdots& \\
    & (t-\theta)^{k-1}\Upsilon_{\phi} &     \\
    (t-\theta)^kf_{\lambda_1}(t) & \cdots & (t-\theta)^kf_{\lambda_r}(t)
\end{pmatrix}\in \Mat_{(rk+1)\times r}(\mathbb{T}).
\]
Indeed, by Lemma \ref{L:coef}, $\mathcal{F}_k$ is equal to $\mathcal{F}$ given in \eqref{D:mathcal F}, namely the matrix whose $i$-th column consists of the entries of the Anderson generating function of $\rho^{[k]}$ with respect to the period $\varpi_{k,i}$ defined in \eqref{E:periodsofext}. Let $\Psi_{\phi}\in \GL_r(\mathbb{T})$ be as in \eqref{E:drinfeldpsi} and consider
\[
\Psi_{k}:=\frac{1}{(\omega_C^{(1)})^k}\Psi_{\phi}\in \GL_r(\mathbb{T}).
\]
Here, again by Lemma \ref{L:coef}, we note that $\Psi_{k}$ is the matrix  $\Psi\in \GL_r(\mathbb{T})$ defined in \cite[\S2.4.1]{GG26}  for $\rho^{[k]}$. 
Finally, we set 
\[
\mathfrak{N}_0:=\Log_{\phi}\frac{1}{\theta^{\mathfrak{m}_{0}}}\in \mathbb{M}_{\phi}.
\]

Our first main result can be seen as a generalization of \cite[Thm. 4.2(i)]{GG26} and it should be compared with both Deligne's construction of $c^\pm(M)$ in \cite{Del} and with Brown's construction of matrices of algebraic Mellin integrals in \cite{Bro26}.
\begin{theorem}\label{T:1} Let $\phi$ be a Drinfeld module defined over $A$ and let $L(\phi^{\vee},k)$ be the dual Goss $L$-function associated to $\phi$ as described in \S \ref{S:3}. Then for each $k\in \mathbb{Z}_{\geq 1}$ ($k= 0$ resp.), there exist $\zz_1,\dots,\zz_r\in \overline{K}^{rk+1}$ ($\zz\in \overline{K}$ resp.) and $\mathfrak{c}_{k},\mathfrak{c}_{0}\in K^{\times} $ such that 
     \begin{equation}\label{E:detlvalue1}
L(\phi^{\vee},k)=\mathfrak{c}_{k}\det\begin{bmatrix}\mathcal{M}_{\zz_1}(\mathfrak{N}_{r(k-1)+2}\cdot \mathcal{F}_k(\Psi_k^{\tr})^{(-1)})&\cdots & \mathcal{M}_{\zz_1}(\mathfrak{N}_{rk+1}\cdot \mathcal{F}_k(\Psi_k^{\tr})^{(-1)})\\
     \vdots & & \vdots \\
\mathcal{M}_{\zz_r}(\mathfrak{N}_{r(k-1)+2}\cdot \mathcal{F}_k(\Psi_k^{\tr})^{(-1)})&\cdots & \mathcal{M}_{\zz_r}(\mathfrak{N}_{rk+1}\cdot \mathcal{F}_k(\Psi_k^{\tr})^{(-1)})
     \end{bmatrix}
     \end{equation}
     and 
     \begin{equation} \label{E:detlvalue2}     L(\phi^{\vee},0)= \mathfrak{c}_{0}\mathcal{M}_{\zz}(\mathfrak{N}_0\cdot \mathcal{F}_0(\Psi_{\phi}^{\tr})^{(-1)}).
     \end{equation} 
     \end{theorem}

The strategy for the proof of Theorem \ref{T:1} can be described as follows. Using \cite[\S5]{Mau22}, we first obtain a set of generators of the period lattice $\Lambda_{\rho^{[k]}}$ as an $\bA$-module so that their last $r$-entries are given in terms of periods and quasi-periods of the Drinfeld module $\phi$ (Proposition \ref{P:periodbasis}). Combining this analysis with \cite[Thm. 1.6]{GG26}, we obtain Theorem \ref{T:mellin} deducing formulas for entries of the logarithms of $\rho^{[k]}$. On the other hand,  in Proposition \ref{P:lvalueTaelman}, we realize  $L(\phi^{\vee},k)$ in terms of the Taelman $L$-value $L(\rho^{[k]}/A)$ of $\rho^{[k]}$ (see \S\ref{S:lfunctionsTaelman} for the definition of Taelman $L$-values). Then, using \cite[Thm. 4.4]{ANDTR20}, we write $L(\rho^{[k]}/A)$ in terms of a determinant of a matrix whose entries consist of last $r$-entries of certain logarithms of $\rho^{[k]}$. Finally, using Theorem \ref{T:mellin}, we conclude the proof of Theorem \ref{T:1}.

Indeed, for certain Drinfeld modules and particular values of $k$, the formulation in Theorem \ref{T:1} is more explicit: Let $\phi$ be a Drinfeld module of rank $r$ as in \eqref{E:Drinfeld} 
so that $k_1,\dots,k_{r-1}\in \mathbb{F}_q$ and $k_r\in \mathbb{F}_q^{\times}$. First, we fix a $(q-1)$-st root of $-\theta$ and define \textit{the Carlitz period} by
\begin{equation}\label{E:tildepi}
\tilde{\pi}:=\theta(-\theta)^{1/(q-1)}\prod_{j=1}^{\infty}\left(1-\theta^{1-q^j}\right)^{-1}\in \mathbb{C}_{\infty}^{\times}.
\end{equation}
Moreover, we continue to let $\lambda_1,\dots,\lambda_r$ be fixed generators of $\Lambda_{\phi}$. For $1\leq \mu \leq r-1$, we set 
\[
\mathcal{P}_{\mu}:=(-1)^{k}\tilde{\pi}^k[F_{\tau^{\mu}}(\lambda_1),\dots,F_{\tau^{\mu}}(\lambda_r)]\in \Mat_{1\times r}(\mathbb{C}_{\infty})
\]
 where $F_{\tau^{\mu}}$ is a quasi-periodic function (see \eqref{D:quasi periodic function} for its explicit definition). We also define 
\[
\mathcal{P}_{r}:=(-1)^k\tilde{\pi}^k[\lambda_1,\dots,\lambda_r] \in \Mat_{1\times r}(\mathbb{C}_{\infty}).
\]
We now state our next result whose proof heavily relies on Theorem \ref{T:mellin}, an analysis of the logarithm function of $\rho^{[k]}$ in \S4 and a generalization of the strategy used to prove \cite[Thm. 1.1]{Gez21}.
\begin{theorem}\label{T:specval} 
     Let $\zz_1,\dots,\zz_r\in \Mat_{(rk+1)\times r}(\mathbb{C}_{\infty})$ ($\zz\in \mathbb{C}_{\infty}$ resp.), $\mathfrak{c}_{k}$ and $\mathfrak{c}_{0}$ be as in Theorem \ref{T:1}, let $\phi$ be a Drinfeld module with the restrictions discussed above and let $L(\phi^{\vee},k)$ be the dual Goss $L$-function associated to $\phi$ as described in \S \ref{S:3}. Let  $0\leq k \leq q-1$. Then 
for $1\leq i \leq r$,   $\zz_i\in \Mat_{(rk+1)\times 1}(\mathbb{F}_q)$ is the $i$-th unit vector and $\zz=\mathfrak{c}_{k}=\mathfrak{c}_{0}=1$. Moreover, for $1\leq j \leq r$, set \[
     \mathcal{L}_j:=\frac{1}{\theta-t}\mathcal{P}_j(\Psi_{k}^{\tr})^{(-1)}\in \Mat_{1\times r}(\mathbb{T}).
     \]    Then, we have 
\[
L(\phi^{\vee},k)=\det\begin{bmatrix}
\mathcal{M}_{\zz_{r(k-1)+2}}(\mathcal{L}_1)&  \cdots & \mathcal{M}_{\zz_{r(k-1)+2}}(\mathcal{L}_r)\\
\vdots&  & \vdots \\
\mathcal{M}_{\zz_{rk+1}}(\mathcal{L}_{1})& \cdots & \mathcal{M}_{\zz_{rk+1}}(\mathcal{L}_r)
	\end{bmatrix}
\]
and 
\[
     L(\phi^{\vee},0)=\mathcal{M}_{1}\left(\frac{1}{\theta-t}[\lambda_1,\dots,\lambda_r](\Psi_{\phi}^{\tr})^{(-1)}\right).    \]
\end{theorem}

\begin{remark}\label{R:KZperiods} In what follows, we briefly describe \textit{periods} which are defined by Kontsevich and Zagier in \cite[\S1.1]{KZ01} (see also \cite{HMS} for more details and alternative equivalent definitions of periods). Let $X$ be a smooth quasi projective variety over $\overline{\mathbb{Q}}$ and $Y$ be a subvariety. Let $\omega$ be a closed algebraic $n$-form on $X$ which is defined over $\overline{\mathbb{Q}}$ vanishing on $Y$. Let $C$ be a singular $n$-chain on $\mathbb{C}$-points $X(\mathbb{C})$ of $X$ whose boundary is contained the set of $\mathbb{C}$-points $Y(\mathbb{C})$ of $Y$. Then we call $\gamma\in \mathbb{C}$ \textit{a period} if
\[
\gamma=\int_{C}\omega.
\]
For a given pure motive $M$ over $\mathbb{Q}$, we know that the entries of the matrix representing the comparison isomorphism between Betti and de Rham cohomology are periods. In other words, Deligne's period $c^{\pm}(M)$ is derived from the determinant of a matrix whose entries consist of periods. In an ongoing work of the authors, motivated by Theorem \ref{T:1}, we aim to define an analogue of Kontsevich-Zagier periods in our setting as an image of the map $\mathcal{M}_{\zz}$ of certain row vectors in the Tate algebra and $\zz$ having entries in $\overline{K}$. Indeed, the matrix $\Psi_{\phi}^{(-1)}$ can be viewed as an element in the Betti cohomology module for $\phi$ (see \cite[\S 1.2]{GP24}). Further, letting $M_{\phi}$ be the $t$-motive of $\phi$, we note that the de Rham cohomology module corresponding to $\phi$ is naturally isomorphic to $\tfrac{1}{\theta-t} M_{\phi}/M_{\phi}$, and in certain cases of our main formula, the quantity $\tfrac{1}{\theta-t}$ naturally appears inside the Mellin transform, as can be seen in the second part of Theorem \ref{T:specval}. We plan to study the exact relationship between our Mellin transform and the pairing between Betti and de Rham cohomology in our upcoming work and use this to give a general definition of Kontsevich-Zagier type periods in the setting of Anderson modules. We hope this will make a more transparent comparison with the classical setting.

\end{remark}

The outline of the present paper can be described as follows. In \S\ref{S:2}, we introduce some fundamental objects, such as Anderson modules, their $t$-motives and $t$-comotives as well as analysis of certain examples of Anderson modules. In \S\ref{S:3}, after introducing Goss $L$-functions of Drinfeld modules, we provide a proof for Theorem \ref{T:1}. In \S\ref{S:4}, we investigate the domain of convergence of $\Log_{\rho^{[k]}}$ and prove Proposition \ref{P:roclog}. In \S\ref{S:4.1} and \S\ref{S:4.2}, after analyzing the norms of the $\bA$-module generators of $\Lambda_{\rho^{[k]}}$ for a certain class of Drinfeld modules with $\mathfrak{m}_0=0$, we finally prove Theorem \ref{T:specval} in \S\ref{S:specvalues}.

\subsection*{Acknowledgments} We would like to thank A. Maurischat for a discussion on his results in \cite{Mau22}. The first author acknowledges support by NSTC Grant 113-2115-M-007-001-MY3. The second author acknowledges support from the NSF under Grant No. 2302399. Data sharing is not applicable for this article because no datasets were generated or analyzed while preparing the present paper.

\section{Preliminaries and Background}\label{S:2}
Our goal in this section is to provide necessary background on the main objects of the present paper. More precisely, we introduce Anderson modules, Anderson generating functions, (abelian) $t$-motives, $t$-comotives as well as the map $\varphi_{\phi}$ constructed for any abelian and uniformizable Anderson module $\phi$. We also provide several examples of Anderson modules and analyze their corresponding $t$-motives and $t$-comotives. For more details, we refer the reader to \cite{And86, BP20, HJ20, NP21}.
\subsection{Anderson modules, Anderson generating functions and hyperderivatives} \label{S:Abelian And Modules} Recall that $\bA=\mathbb{F}_q[t]$. We call $L$ is \textit{an $\bA$-field} if there exists a ring homomorphism $\ch: \bA\to L$.
\begin{definition}
   \begin{itemize} 
       \item[(i)] \textit{An Anderson module of dimension $d$} (defined over an $\bA$-field $L$) is an $\mathbb{F}_q$-algebra homomorphism $\phi:\bA\to \Mat_d(L[\tau])$ given by 
       \[
       \phi_{t}=\dmat_{\phi}[t]+A_1\tau+\cdots+A_{s}\tau^s
       \]
       for some $s\in \ZZ_{\geq 0}$ so that $\dmat_{\phi}[t]:=\ch(t) \Id_d+N$ where $N\in \Mat_d(L)$ is a nilpotent matrix. 
       \item[(ii)] Let $\psi$ be an Anderson module of dimension $d'$.  The morphisms between $\phi$ and $\psi$ are given by matrices $U\in \Mat_{d\times d'}(L[\tau])$ satisfying
       \[
       \phi_{t}U=U\psi_{t}.
       \]
        When such a $U$ exists, we say that $\phi$ and $\psi$  are \textit{isogenous} (over $L$). Furthermore, if $U\in \GL_d(L[\tau])$, then  we say that $\phi$ and $\psi$ are \textit{isomorphic} (over $L$).   \end{itemize} 
\end{definition}
We let $\Lie(\phi)(L)$ be the $\bA$-module $L^d$ whose $\bA$-module structure is given by
\[
a\cdot x:=\dmat_{\phi}[a]x, \ \ a\in \bA, \ \ \ \ x\in L^d.
\]
For $1\leq i \leq d$, we say that the $i$-th coordinate of $\phi$ is \textit{tractable} if, for any $x\in \mathbb{C}_{\infty}^d$ and $a\in \bA$, the $i$-th coordinate of $a\cdot x\in \Lie(\phi)(\mathbb{C}_{\infty})$ is equal to the $i$-th coordinate of $(a|_{t=\theta})x$.
We also let $\phi(L)$ be the $\bA$-module $L^d$ whose $\bA$-module structure is given by
\[
a\cdot x:=\phi_a x, \ \ a\in \bA, \ \  \  x\in L^d.
\]

Until the end of this subsection, we assume that $L$ is a subfield of $\mathbb{C}_{\infty}$ containing $K$. \textit{The exponential series of $\phi$} is given by
the unique infinite series $\Exp_{\phi}=\sum_{i\geq 0}\alpha_i\tau^i\in \Mat_d(L[[\tau]])$ such that $\alpha_0=\Id_d$ and it satisfies 
\[
\Exp_{\phi}\dmat_{\phi}[t]=\phi_{t}\Exp_{\phi}
\]
in $\Mat_d(L[[\tau]])$. Note that it induces an everywhere convergent $\mathbb{F}_q$-linear function $\Exp_{\phi}:\Lie(\mathbb{C}_{\infty})\to \phi(\mathbb{C}_{\infty})$ defined by
\[
\Exp_{\phi}(\zz)=\Exp_{\phi}\cdot \zz=\sum_{i\geq 0}\alpha_i\zz^{(i)},\ \ \zz\in \Lie(\phi)(\mathbb{C}_{\infty}).
\]
If $\Exp_{\phi}$ is surjective, then we say that $\phi$ is \textit{uniformizable} (see  \cite[\S2.2]{And86} for an example of Anderson module $\phi$ so that $\Exp_{\phi}$ is not surjective). We let $\Lambda_{\phi}\subset\Lie(\phi)(\mathbb{C}_{\infty}) $ be the kernel of $\Exp_{\phi}$ and call $\Lambda_{\phi}$ \textit{the period lattice of $\phi$}. Note that it is a free, finitely generated and discrete $\bA$-module. We also call its non-zero elements  \textit{periods of $\phi$}. 
We further extend the domain of $\Exp_{\phi}$ and, by a slight abuse of notation, let $\Exp_{\phi}: \mathbb{T}^d\to \mathbb{T}^{d}$ so that 
\[
\Exp_{\phi}(f)=\sum_{i\geq 0}\alpha_if^{(i)},\ \ f\in \mathbb{T}.
\]

We define \textit{the logarithm series $\Log_{\phi}=\sum_{i\geq 0} \gamma_i \tau^i\in \Mat_d(L[[\tau]])$ of $\phi$} to be the formal inverse of $\Exp_{\phi}$. It satisfies the functional equation 
\begin{equation}\label{E:logfunc}
    \Log_{\phi}\phi_t=\dmat_{\phi}[t]\Log_{\phi}.
\end{equation}
Note that, there exists a domain $\mathcal{D}\subset \phi(\mathbb{C}_{\infty})$ such that the logarithm series induces an $\mathbb{F}_q$-linear function $\Log_{\phi}:\mathcal{D}\to \Lie(\phi)(\mathbb{C}_{\infty})$ which is given by
\[
\Log_{\phi}(\zz)=\sum_{i\geq 0}\beta_i\zz^{(i)},\ \ \zz\in \mathcal{D}\subset \phi(\mathbb{C}_{\infty}).
\]
We also extend the domain of $\Log_{\phi}$ and consider $\Log_{\phi}: \mathcal{D}'\to \mathbb{T}^{d}$ so that 
\begin{equation}\label{E:extension}
\Log_{\phi}(f)=\sum_{i\geq 0}\beta_if^{(i)},\ \ f\in \mathcal{D}'
\end{equation}
where $\mathfrak{D}'$ is a subdomain of $\mathbb{T}^d$ containing $\mathcal{D}$.

Let $\zz\in \mathbb{C}_{\infty}^d$. \textit{The Anderson generating function $\mathcal{Y}_{\zz}(t)$ of $\phi$ with respect to $\zz$} is the infinite sum given by 
\[
\mathcal{Y}_{\zz}(t):=\sum_{i=0}^{\infty}\Exp_{\phi}\left(\dmat_{\phi}[t]^{-i-1}\zz\right)t^i\in \mathbb{C}_{\infty}[[t]]^d.
\]
Let $P=P_0+P_1\tau+\dots+P_{\ell}\tau^{\ell}\in \Mat_{m\times \ell}(\mathbb{C}_{\infty}[\tau])$ and $\mathcal{G}\in \Mat_{\ell \times v}(\mathbb{T})$. Following the notation in \cite[(2.1.2)]{NP21}, we set 
\[
\langle U|\mathcal{G}\rangle:= U_0\mathcal{G}+B_1\mathcal{G}^{(1)}+\cdots+U_{\ell}\mathcal{G}^{(\ell)}\in \Mat_{m\times v}(\mathbb{T}).
\]
The following properties of Anderson generating functions will be useful throughout the present paper.
\begin{lemma}\cite[Prop. 4.2.7, Lem. 4.2.2, Prop. 4.2.12]{NP21} \label{L:AGF} The following statements hold.
\begin{itemize}
    \item[(i)] We have the identity
    \[
    \mathcal{Y}_{\zz}(t)=\sum_{i=0}^{\infty}\alpha_i\left((\dmat_{\phi}[t]-t\Id_d)^{-1}\right)^{(i)}\zz^{(i)}
    \]
    which converges in $\mathbb{T}^d$.
    \item[(ii)] Let $\mathcal{Y}_{\zz}(t)=((\mathcal{Y}_{\zz})_1,\dots,(\mathcal{Y}_{\zz})_d)^{\tr}\in \mathbb{T}^d.$ Then
    \[
    \Res_{t=\theta}\mathcal{Y}_{\zz}:=(\Res_{t=\theta}(\mathcal{Y}_{\zz})_1dt,\dots, \Res_{t=\theta}(\mathcal{Y}_{\zz})_ddt)^{\tr}=-\zz.
    \]
    \item[(iii)] Let $\lambda \in \Lambda_{\phi}$. Then
    \[
    \langle\phi_{t}|\mathcal{Y}_{\lambda}(t)\rangle=t\mathcal{Y}_{\lambda}(t).
    \]
\end{itemize}
\end{lemma}

We finish this subsection with a brief discussion on hyperderivatives. For more details, we refer the reader to \cite[\S2.4]{NP21} and \cite[\S2]{Pap15}.

Consider the field of formal Laurent series $\CC_{\infty}((t))$  in $t$ with coefficients in $\CC_{\infty}$. Let $n$ be a non-negative integer. \textit{The $n$-th hyperdifferential operator $\partial^{j}_{t}:\CC_{\infty}((t))\to \CC_{\infty}((t))$ with respect to $t$} is given by
	\[
	\partial^{j}_{t}\bigg(\sum_{i=i_0}^{\infty}c_it^i\bigg):=\sum_{i=i_0}^{\infty}\binom{i}{j}c_it^{i-j}\text{, }c_i\in \CC_{\infty}.
	\]
	When $n=0$, we have $\partial_{t}^{0}(f)=f$ for any $f\in \CC_{\infty}((t))$. Let $f_1,\dots,f_{\ell}\in \mathbb{C}_{\infty}((t))$. Then the Leibniz rule for hyperderivatives can be described as
	\begin{equation}\label{E:Leibniz}
	\partial_{t}^{n}(f_1\cdots f_{\ell})=\sum_{\substack{n_1,\dots,n_{\ell}\geq 0\\n_1+\cdots+n_{\ell}=n}}\partial_{t}^{n_1}(f_1)\cdots \partial_{t}^{n_{\ell}}(f_{\ell}).
	\end{equation}

The following lemma can be obtained by using the definition of hyperderivatives and we leave its details to the reader. Before stating it, we note that, when $n=-m$ is a negative integer for $m>0$, we set
\[
\binom{n}{j}:=(-1)^j\binom{m+j-1}{j}.
\]
\begin{lemma}[{see \cite[Lem. 2.3.23, Cor. 2.4.7]{Pap15}}] \label{L:hyp} The following statements hold.
\begin{itemize}
    \item[(i)] Let $f\in \mathbb{F}_q(t)$ and $r\geq 1$ be such that 
    \[
    \partial_t^j(f)=\begin{cases}
        c\in \mathbb{F}_q & \text{ if } j=1\\
        0 & \text{ if } 2\leq j \leq r.
    \end{cases}
    \]
    Then for any $1\leq j \leq r$ and $m\in \mathbb{Z}$, we have 
    \[
    \partial_t^j(f^m)=c^j\binom{m}{j}f^{m-j}.
    \]
    \item[(ii)] For any $i,j\geq 1$, we have
    \[
    \partial_t^j\left(\frac{1}{(t-\theta^{p^i})^m}\right)=(-1)^j\binom{m+j-1}{j}\frac{1}{(t-\theta^{p^i})^{m+j}}.   \]
\end{itemize}
    
\end{lemma}
    
\subsection{\texorpdfstring{$t$}{t}-motives } \label{S:tmotive} 
Our goal in this subsection is to introduce $t$-motives as well as their relation with Anderson modules. We refer the reader to  \cite[\S2.2,\S2.3]{GN25}, \cite{Taelman2} and \cite{Taelman3} for further details on the subject.
	
Let $L$ be an $A$-algebra. For any $f=\sum_{i\geq 0}a_it^i\in L[t]$ and $j\in \ZZ$, we set 
\[
f^{(j)}:=\sum_{i\geq 0}a_i^{q^j}t^i\in L[t].
\]
Furthermore, we denote by $L[t,\tau]:=L[t][\tau]$ the non-commutative ring defined subject to the condition
	\[
	\tau f =f^{(1)}\tau ,\ \ \ \  f\in L[t].
	\]	
	\begin{definition}
		\begin{itemize}
			\item[(i)] \textit{A $t$-motive $M$ over $L$} is a left  $L[t,\tau]$-module which is free and finitely generated over $L[t]$ and $L[\tau]$ so that  there exists a non-negative integer $\nu$ satisfying
		\[
		(t-\theta)^{\nu}M\subset \tau M.
		\] 
            \item[(ii)] The morphisms between  $t$-motives are given by left $L[t,\tau]$-module homomorphisms.  
			\item[(iii)] The tensor product $M_1\otimes_{L[t]} M_2$ of  $t$-motives $M_1$ and $M_2$ is a $t$-motive on which $\tau$ acts diagonally.
			
		\end{itemize}
	\end{definition}

For the rest of this subsection, we let $L$ be a subfield of $\mathbb{C}_{\infty}$ containing $K$.	We briefly describe how one can form an Anderson module from a $t$-motive $M$ over $L$. Let $\{n_1,\dots,n_d\}$ be an $L[\tau]$-basis for $M$. Then, there exists a matrix $\mathcal{N}\in \Mat_d(L)[\tau]$ so that 
	\[
	t\cdot \begin{bmatrix}
	n_1\\\vdots \\ n_d
	\end{bmatrix}=\mathcal{N}\begin{bmatrix}
	n_1\\\vdots \\ n_d
	\end{bmatrix}.
	\]

If we set $\phi:\bA\to \Mat_d(L)[\tau]$ to be the $\mathbb{F}_q$-algebra homomorphism given by $\phi_{t}=\mathcal{N}$ then we see that $\phi$ forms an  Anderson module. An Anderson module constructed via a $t$-motive is called \textit{abelian}. On the other hand, consider $M_{\phi}:=\Mat_{1\times d}(L[\tau])$. It is naturally a free $L[\tau]$-module. We equip $M_{\phi}$ with an $L[t]$-module structure given by 
\begin{equation}\label{E:tmodstr}
\alpha t^i\cdot m:=\alpha  m  \phi_{\theta^i}, \ \ m\in M_\phi, \ \ \alpha\in L.
\end{equation}
Then $M_{\phi}$ forms a left  $L[t,\tau]$-module. We call $M_{\phi}$ \textit{the $t$-motive corresponding to $\phi$}. One can see that $M_{\phi}$ forms a left  $L[t,\tau]$-module which is free and finitely generated over $L[t]$ and $L[\tau]$ so that  
		\[
		(t-\theta)^{d}M\subset \tau M.
		\] 
In other words, $M_{\phi}$ is  a $t$-motive. This correspondence forms an anti-equivalence between the category of $t$-motives over $L$ and the category of abelian Anderson modules defined over $L$  (see \cite[Thm. 1]{Taelman3} and \cite[Thm. 1]{And86}).
        
Let $\{m_1,\dots,m_r\}$ be an $L[t]$-basis for $M$ and let $\Theta\in  \Mat_r(L[t])$ be such that 
\begin{equation}\label{E:003}
\tau \cdot \begin{bmatrix}m_1\\ \vdots\\
m_r
\end{bmatrix}=\Theta\begin{bmatrix}m_1\\ \vdots\\
m_r
\end{bmatrix}.
\end{equation}
We call  $M$ \textit{rigid analytically trivial} if there exists $\Upsilon\in \GL_r(\mathbb{T})$ satisfying
\[
\Upsilon^{(1)}=\Theta\Upsilon.
\]
We further call $\Upsilon$ \textit{a rigid analytic trivialization of $M$}.

\subsection{$t$-comotives} \label{S:dualmot}  We now assume that $K\subset L\subseteq \mathbb{C}_{\infty}$ is a perfect field and define the non-commutative polynomial ring $\Mat_{d}(L[\sigma])$ subject to the relation
\[
\sigma B=B^{(-1)} \sigma, \ \ B\in \Mat_{d}(L).
\]
We set $L[t,\sigma]$ to be the ring of polynomials in $t$ and $\sigma$ with coefficients in $L$ subject to the condition
\[
tx=xt, \ \ t\sigma=\sigma t, \ \ \sigma x=x^{(-1)}\sigma=x^{1/q} \sigma, \ \ x\in L.
\]
For any $ g= \sum_{i\geq 0}c_i\tau^i\in L[\tau]$, we also set 
\[
g^{*}:=\sum_{i\geq 0}c_i^{(-i)}\sigma^i\in L[\sigma].
\]
For any $B=(b_{ij})_{ij}\in \Mat_d(L[\tau])$, we extend the above definition to matrices so that
\[
B^{*}:=(b_{ij}^{*})^{\tr}_{ij}\in \Mat_d(L[\sigma]).
\]
\begin{definition}
	\begin{itemize}
		\item[(i)] \textit{A ($t$-finite) $t$-comotive $N$} is a left $L[t,\sigma]$-module which is free and finitely generated over $L[t]$ and $L[\sigma]$ such that there exists $\mu\in \mathbb{Z}_{\geq 0}$ satisfying
		\[
		(t-\theta)^{\mu} N\subset \sigma N.
		\] 
		\item[(ii)] Let $N_1$ and $N_2$ be $t$-comotives. Morphisms between $N_1$ and $N_2$ are given by left $L[t,\sigma]$-module homomorphisms.
		\item[(iii)] The tensor product of $N_1$ and $N_2$ forms a $t$-comotive  $N_1\otimes N_2:=N_1\otimes_{L[t]} N_2$ where $\sigma$ acts diagonally.
	\end{itemize}
\end{definition}

Let $\{\bn_1,\dots,\bn_r\}$ be an $L[t]$-basis for $N$ and let $\Phi\in  \Mat_r(L[t])$ be such that 
\[
\sigma \cdot \begin{bmatrix}\bn_1\\ \vdots\\
\bn_r
\end{bmatrix}=\Phi\begin{bmatrix}\bn_1\\ \vdots\\
\bn_r
\end{bmatrix}.
\] 
We say that $N$ is \textit{rigid analytically trivial} if there exists $\Psi\in \GL_r(\mathbb{T})$ satisfying 
\[
\Psi^{(-1)}=\Phi\Psi.
\]
We call $\Psi$ \textit{a rigid analytic trivialization of $N$}.

Anderson, in his unpublished notes, constructed an equivalence between a subcategory of Anderson modules, that is the category of \textit{$t$-finite Anderson modules} and the category of $t$-comotives. We refer the reader to \cite[\S2.5.2]{HJ20} for more details of this correspondence. Let $\phi$ be an Anderson module of dimension $d$ which is defined over $L$. We set $N_{\phi}:=\Mat_{1\times d}(L[\sigma])$ and call $N_{\phi}$ \textit{the $t$-comotive of $\phi$}. It is naturally a free and finitely generated $L[\sigma]$-module. It is also equipped with an  $L[t]$-module structure given by 
\[
\alpha t^i\cdot n:=\alpha  n \phi_{\theta^i}^{*}, \ \ n\in N_\phi, \ \ \alpha\in L.
\]
Hence, $N_{\phi}$ forms a left $L[t,\sigma]$-module. If $N_{\phi}$ is also a free $L[t]$-module of finite rank, we say that $\phi$ is \textit{a $t$-finite Anderson module}. It is important to note that, by \cite[Thm. A]{Mau21},  being an abelian Anderson module is equivalent to being $t$-finite.

\subsection{The map $\varphi_{\phi}$}\label{S:map phi}
In this subsection we give the details necessary for the definition of the map $\varphi_{\phi}$ introduced in \eqref{E:var phi intro}. The details of the map $\varphi_{\phi}$ were developed in \cite{GG25} and \cite{GG26}, and so we only give a brief overview here, and refer the reader to loc. cit. for many of the proofs and details. We maintain the notation and definitions of the previous two subsections.

Let $M$ be the $t$-motive of an abelian Anderson module $\phi$. Let $\bm:=[\bm_1,\dots,\bm_r]^{\tr}\in \Mat_{r\times 1}(M)$ be a certain $\mathbb{C}_{\infty}[t]$-basis for $M$ (see \cite[2.4.1]{GG26}) so that 
\[
\tau \cdot \bm=\Phi \bm
\]
and let $\mathbf{b}:=\{\bb_1,\dots,\bb_d\}$ be an $L[\tau]$-basis for $M$. Then for each $1\leq \mu \leq r$ and $1\leq \ell \leq d$, there exist  elements $ \mathcal{C}_{\ell}^{[\mu]}\in \mathbb{C}_{\infty}[\tau] $ such that 
\begin{equation}\label{E:004}
\mathbf{m}_{\mu}=\sum_{\ell=1}^{d}\mathcal{C}_{\ell}^{[\mu]}\bb_{\ell}
\end{equation}
and for each $1\leq \ell \leq d$ and $1\leq \mu \leq r$, there exist elements 
$
\mathcal{B}_{\mu}^{[\ell]}\in \mathbb{C}_{\infty}[t]
$ such that
\begin{equation}\label{E:008}
   \bb_{\ell}=\sum_{\mu=1}^{r}\mathcal{B}_{\mu}^{[\ell]}\mathbf{m}_{\mu}.
\end{equation}

We remind the reader the definitions of $\mathbb{M}_{\phi}$ and its norm from \eqref{D:M phi}. Set $\mathfrak{S}_{-1}:=\Id_r$ and for $n\geq 0$, set $\mathfrak{S}_{n}:=(\Phi^{\tr})^{(n)}\cdots (\Phi^{\tr})^{(1)}\Phi^{\tr} \in \GL_r(\mathbb{C}_{\infty}(t))\cap\GL_r(\mathbb{T}) $. For each $1\leq \mu \leq  r$, let $\mathfrak{e}_{\mu}\in \Mat_{r\times 1}(\mathbb{F}_q)$ be the $\mu$-th unit vector. Consider the map
$\varphi_{\phi}:\mathbb{M}_{\phi}\to \Mat_{1\times r}(\mathbb{T})$ given by   
\begin{equation}\label{D:varphi def}
\varphi_{\phi}\left(\left(\sum_{n=0}^{\infty}a_{1,n}\tau^n,\dots,\sum_{n=0}^{\infty}a_{d,n}\tau^n\right)\right):=\sum_{\ell=1}^d\sum_{n=0}^{\infty}a_{\ell,n}\sum_{\mu=1}^{r}(\mathcal{B}_{\mu}^{[\ell]})^{(n)}\mathfrak{e}_{\mu}^{\tr}\mathfrak{S}_{n-1}.
\end{equation}
In \cite[Thm. 3.7]{GG25}, we show that $\varphi_\phi$ is well defined, continuous and injective.
\begin{remark}
We give a brief remark to clarify for the reader the level of generality allowed when selecting the bases in this section, as well to explain how our main results depend on the choice of bases. The $\C_\infty[t]$-basis $\bm$ and $\C_\infty[\tau]$-basis $\bb$ for a $t$-motive $M$ may be chosen arbitrarily. Then, the choice of another basis $\tilde{\bm}$ is prescribed. Notably, this fixes the structure matrix $\Phi$, which appears in the definition of $\mathfrak{S}_n$, which is used in the definition of $\varphi_{\phi}$. The values $\mathcal{B}_\mu^{[\ell]}$ also depend on these choices. Ultimately, our construction gives a formula for the logarithm function $\Log_{\phi}$, and is independent of the choice of $\C_\infty[t]$-bases $\bm$ (and $\tilde{\bm}$). On the other hand, choosing a different $\C_\infty[\tau]$-basis will yield $\Log_{\phi}$ for an isomorphic Anderson module. We refer the reader to \cite[Rem. 3.10]{Gre24} for the exact relationship between these isomorphic Anderson modules.
\end{remark}

\subsection{The case of Drinfeld modules}  \label{S:Drinfeld} Let $\phi$ be a Drinfeld module of rank $r$ defined over $L$ as in \eqref{E:Drinfeld}. Consider the $t$-motive $M_{\phi}$ given as in \S\ref{S:tmotive}. We define the matrix
\begin{equation}\label{E:phitheta}
\Theta_{\phi}:=\begin{pmatrix}
&1& & &  \\
& & \ddots & & \\
& & & \ddots & \\
& & & &  1\\
\frac{t-\theta}{k_r}&-\frac{k_1}{k_r} & \dots & \dots & -\frac{k_{r-1}}{k_r}
\end{pmatrix}\in \GL_r(\mathbb{T})\cap \Mat_r(L[t]).
\end{equation}
One can choose an $L[t]$-basis $\bm=[\bm_1,\dots,\bm_{r}]^{\tr}\in \Mat_{r\times 1}(M_{\phi}) $ so that $\bm_i:=\tau^{i-1}$ for $1\leq i \leq r$. We see that 
\[
\tau \cdot \bm=\Theta_{\phi} \bm.
\]
Note that $\phi$ is always uniformizable and its period lattice $\Lambda_{\phi}$ is a free and finitely generated $\bA$-module of rank $r$. We let $\{\lambda_1,\dots,\lambda_r\}$ be a set of generators of $\Lambda_{\phi}$ as an $\bA$-module. For any $i\in \{1,\dots,r\}$, we denote by $f_{\lambda_i}(t)$  the Anderson generating function of $\phi$ with respect to $\lambda_i$. We define 
\begin{equation}\label{E:Upsilon}
\Upsilon_{\phi}:=\begin{pmatrix} f_{\lambda_1}(t) & \cdots & \cdots & f_{\lambda_r}(t)\\
f_{\lambda_1}(t)^{(1)} & \cdots & \cdots & f_{\lambda_r}(t)^{(1)}\\
\vdots  & & & \vdots\\
f_{\lambda_1}(t)^{(r-1)} & \cdots & \cdots &f_{\lambda_r}(t)^{(r-1)}
\end{pmatrix}\in \Mat_{r\times r}(\mathbb{T}).
\end{equation}
By \cite[\S4.2]{Pel08} (see also Lemma \ref{L:AGF}(iii)), we know that $\Upsilon_{\phi}\in \GL_r(\mathbb{T})$ and moreover it satisfies 
\begin{equation}\label{E:RAT1}
\Upsilon_{\phi}^{(1)}=\Theta_{\phi}\Upsilon_{\phi}.
\end{equation}
Hence $M_{\phi}$ is rigid analytically trivial. 

Next, we define
\begin{equation}\label{E:V def}
V:=\begin{pmatrix}
k_1&k_2^{(-1)} &k_3^{(-2)}& \dots &k_r^{(1-r)}\\
\vdots &\vdots&\vdots&  \iddots & \\
\vdots & \vdots& k_r^{(-2)} & & \\
\vdots &k_r^{(-1)} & & & \\
k_r &  & & & 
\end{pmatrix}\in \GL_r(L).
\end{equation}
Observe that  \begin{equation}\label{E:Phi def Drinfeld modules}
\Phi_{\phi}:=(V^{(-1)})^{-1}\Theta_{\phi}^{\tr}V=\begin{pmatrix}
&1& & &  \\
& & \ddots & & \\
& & & \ddots & \\
& & & &  1\\
\frac{t-\theta}{k_r^{(-r)}}&-\frac{k_1^{(-1)}}{k_r^{(-r)}} & \dots & \dots & -\frac{k_{r-1}^{(-(r-1))}}{k_r^{(-r)}}
\end{pmatrix}\in  \GL_r(\mathbb{T})\cap \Mat_r(L[t]).
\end{equation}
We now assume that $L$ is a perfect field and recall the $t$-comotive $N_{\phi}$ of $\phi$ defined in \S\ref{S:dualmot}. We can let $\bn=[1,\dots,\sigma^{r-1}]^{\tr}\in \Mat_{r\times 1}(N_{\phi})$ so that 
\[
\sigma \cdot \bn=\Phi_{\phi} \bn.
\]
We see that $\bn$ forms an $L[t]$-basis for $N_{\phi}$. Moreover, $\{1\}$ forms an $L[\sigma]$-basis for $N_{\phi}$. 
Consider  
\begin{equation}\label{E:drinfeldpsi}
\Psi_{\phi}:=V^{-1}((\Upsilon_{\phi}^{(1)})^{\tr})^{-1}\in \GL_r(\mathbb{T}).
\end{equation}
Thus, we obtain
\begin{equation}\label{E:RAT2}
\Psi_{\phi}^{(-1)}=\Phi_{\phi}\Psi_{\phi}.
\end{equation}
In other words, $\Psi_{\phi}$ is a rigid analytic trivialization of $N_{\phi}$ and hence $N_{\phi}$ is rigid analytically trivial. 

Finally, we describe the map $\varphi_{\phi}$ in the Drinfeld module case. Note that in this case
\[
\mathbb{M}_{\phi}:=\left\{\sum_{n=0}^{\infty}a_n\tau^n \ \ | \ \ \text{ for }1\leq j \leq r, \ \  |a_n|||f_{\lambda_j}(t)^{(n)}||\to 0 \text{ as } n\to \infty\right\}.
\]
We also note that $\mathfrak{S}_{-1}=\Id_r$ and for $n\geq 0$, set $\mathfrak{S}_{n}=(\Phi_{\phi}^{\tr})^{(n)}\cdots (\Phi_{\phi}^{\tr})^{(1)}\Phi_{\phi}^{\tr} \in \GL_r(\mathbb{C}_{\infty}(t))\cap\GL_r(\mathbb{T}) $. The map
$\varphi_{\phi}:\mathbb{M}_{\phi}\to \Mat_{1\times r}(\mathbb{T})$ is given by   
\[
\varphi_{\phi}\left(\sum_{n=0}^{\infty}a_{n}\tau^n\right)=\sum_{n=0}^{\infty}a_{n}((k_r^{(-r)})^{-1})^{(n)}\mathfrak{e}_{r}^{\tr}\mathfrak{S}_{n-1}.
\]

\subsection{Carlitz tensor powers}\label{S:Carlitztensorcase} We  call the Drinfeld module $C$ of rank one given by $C_{t}:=\theta+\tau$ \textit{the Carlitz module}. Let $k\in \mathbb{Z}_{\geq 1}$. We consider the left $L[t,\tau]$-module 
\[
M_{C^{\otimes k}}:=M_{C}\otimes_{L[t]} \cdots \otimes_{L[t]} M_{C}=L[\tau]\otimes_{L[t]}\cdots \otimes_{L[t]} L[\tau]
\]
so that $\tau$ acts diagonally. Let $\bm_1$ be an  $L[t]$-basis for $M_C$. Then $\bm:=\bm_1\otimes \cdots \otimes \bm_1\in M_{C^{\otimes k}}$ is an $L[t]$-basis for $M_{C^{\otimes k}}$ so that 
\[
\tau \cdot \bm=(t-\theta)^k \bm.
\]
Furthermore, the set $\{\bm, (t-\theta)\bm,\dots,(t-\theta)^{k-1}\bm\}$ forms an $L[\tau]$-basis for $M_{C^{\otimes k}}$ and hence it is of dimension $k$ over $L[\tau]$. Therefore, $M_{C^{\otimes k}}\cong \Mat_{1\times k}(L[\tau])$ as $L[\tau]$-modules. The corresponding Anderson module $C^{\otimes k}:\bA\to \Mat_k(L[\tau])$ is given by
\[
C^{\otimes k}_{t}:=\begin{pmatrix}
		\theta&1& & \\
		& \ddots&\ddots & \\
		& & \ddots & 1\\
		& & & \theta  
		\end{pmatrix}+\begin{pmatrix}
		& &  & &\\
		& & & &\\
		& &  & &\\
		1&&&  &
		\end{pmatrix}\tau.
		\]
   Since $M_{C^{\otimes k}}$ is of rank one over $L[t]$, $C^{\otimes k}$ is an abelian Anderson module.
   
We fix a $(q-1)$-st root of $-\theta$. \textit{The Anderson-Thakur function} $\omega_C$ is the Anderson generating function of $C$ with respect to the Carlitz period $\tilde{\pi}$ defined in \eqref{E:tildepi}. It has an infinite product expansion given by 
\begin{equation}\label{D:omega_C}
\omega_C:=(-\theta)^{1/(q-1)}\prod_{j=0}^{\infty}\left(1-\frac{t}{\theta^{q^j}}\right)^{-1}\in \mathbb{T}^{\times}.
\end{equation}
We let $\mathbb{S}_{-1}:=1$ and for $n\geq 0$, define  
\begin{equation}\label{E:sn}
\mathbb{S}_n:=(t-\theta^{q^n})\cdots (t-\theta)\in A[t].
\end{equation}
For $\ell\geq 1$, Lemma \ref{L:AGF}(iii) implies
\begin{equation}\label{E:AGF4}
(\omega_C^k)^{(\ell)}=\mathbb{S}_{\ell-1}^k\omega_C^k.
\end{equation}
Thus, if we let $\ell=1$ in \eqref{E:AGF4}, we see that $\omega_C^k$ is a rigid analytic trivialization for $M_{C^{\otimes k}}$, implying that $M_{C^{\otimes k}}$ is rigid analytically trivial. Hence, by \cite[Thm. 4]{And86}, $C^{\otimes k}$ is uniformizable and its period lattice $\Lambda_{C^{\otimes k}}$ is a rank one $\bA$-module.

We assume that $L$ is a perfect field and consider the $L[t,\sigma]$-module given by
\[
N_{C^{\otimes k}}:=N_{C}\otimes_{L[t]} \cdots \otimes_{L[t]} N_{C}=L[\sigma]\otimes_{L[t]}\cdots \otimes_{L[t]} L[\sigma]
\]
where $\sigma$ acts diagonally.  We let $\bn_1$ be an $L[t]$-basis for $N_C$. Then $\bn_k=\bn_1\otimes \cdots \otimes \bn_1\in N_{C^{\otimes k}}$ is an $L[t]$-basis for $N_{C^{\otimes k}}$ so that 
\[
\tau \bn_k=(t-\theta)^k \bn_k.
\]
Moreover, the set $\{\bn_k, (t-\theta)\bn_k,\dots,(t-\theta)^{k-1}\bn_k\}$ forms an $L[\tau]$-basis for $N_{C^{\otimes k}}$ and hence it is of dimension $k$ over $L[\sigma]$. Therefore, $N_{C^{\otimes k}}\cong \Mat_{1\times k}(L[\sigma])$ as $L[\sigma]$-modules.
We now define the element 
$
\Omega(t):=(\omega_C^{(1)})^{-1}\in \mathbb{T}^{\times}.
$
Note that 
\begin{equation}\label{E:RAT3}
(\Omega(t)^k)^{(-1)}=(t-\theta)^k\Omega(t)
\end{equation}
implying that $\Omega(t)^k$ is a rigid analytic trivialization of $N_{C^{\otimes k}}$. That is, $N_{C^{\otimes k}}$ is rigid analytically trivial.

We now describe the map $\varphi_{C^{\otimes k}}$ in the Carlitz tensor power case. Let $\mathfrak{Y}(t)=(\mathfrak{Y}_1(t),\dots,\mathfrak{Y}_k(t))^{\tr}$ be the Anderson generating function of $C^{\otimes k}$ with respect to a fixed generator of $\Lambda_{C^{\otimes k}}$. Note that 
\[
\mathbb{M}_{C^{\otimes k}}:=\left\{\left( \sum_{n=0}^{\infty}a_{1,n}\tau^n, \dots, \sum_{n=0}^{\infty}a_{k,n}\tau^n \right)  \ \ | \ \ \text{ for } 1\leq \ell \leq k, \ \ |a_{\ell,n}| ||\mathfrak{Y}_{\ell}(t)^{(n)} ||\to 0 \text{ as } n\to \infty\right\}.
\]
We also note that $\mathfrak{S}_{-1}=\Id_r$ and for $n\geq 0$, set $\mathfrak{S}_{n}=(t-\theta^{q^n})\cdots (t-\theta) \in \mathbb{C}_{\infty}[t]\cap \mathbb{T}$.  The map
$\varphi_{C^{\otimes k}}:\mathbb{M}_{C^{\otimes k}}\to \mathbb{T}$ is given by   
\[
\varphi_{\phi}\left(\left( \sum_{n=0}^{\infty}a_{1,n}\tau^n, \dots, \sum_{n=0}^{\infty}a_{k,n}\tau^n \right)\right)=\sum_{\ell=1}^k\sum_{n=0}^{\infty}a_{\ell,n}(t-\theta^{q^n})^{\ell-1}\mathfrak{S}_{n-1}.
\]
\subsection{The tensor product of Drinfeld modules with Carlitz tensor powers}
\label{S:drinfeldtens} Let $\phi$ be a Drinfeld module of rank $r\geq 1$ defined as in \eqref{E:Drinfeld}. By \S\ref{S:tmotive}, we know that $M_{\phi}$ is a free and finitely generated $L[t]$-module of rank $r$. 
We let $M_0=M_{\phi}$ and for $k\geq 1$, set  $M_k:=M_{\phi}\otimes_{L[t]}M_{C^{\otimes k}}$.  Then $\mathfrak{F}_k:=[\mathfrak{m}_1\otimes \bm,\dots,\mathfrak{m}_r\otimes \bm]^{\tr}\in \Mat_{r\times 1}(M_k)$  is an $L[t]$-basis for $M_k$ so that 
\begin{equation}\label{E:extmapk}
\tau \cdot \mathfrak{F}_k=(t-\theta)^k\Theta_{\phi}\mathfrak{F}_k.
\end{equation}
Moreover, there exists an $L[\tau]$-basis $\bb_k\in \Mat_{(rk+1)\times 1}(M_k)$ of $M_k$ which is given by 
        \begin{multline*}
       \bb_k:=[\mathfrak{m}_1\otimes \bm, \dots, \mathfrak{m}_r\otimes \bm
		, (t-\theta)\mathfrak{m}_1\otimes \bm, \dots,(t-\theta) \mathfrak{m}_r\otimes \bm,\dots,\\
         (t-\theta)^{k-1}\mathfrak{m}_1\otimes \bm, \dots, (t-\theta)^{k-1}\mathfrak{m}_r\otimes \bm_1, (t-\theta)^{k}\mathfrak{m}_1\otimes \bm]^{\tr}.
        \end{multline*}
To ease the notation, we denote by $\rho^{[k]}$ the abelian Anderson module corresponding to the $t$-motive $M_k$. More precisely, the Anderson module 
		\[
		\rho^{[k]}:\bA\to \Mat_{rk+1}(\mathbb{C}_{\infty})[\tau]
		\] is given by
		\[
		\rho^{[k]}_{t}:=	\begin{pmatrix}
		\theta&\cdots &0&1&0& \cdots & 0  \\
		& \ddots& & \ddots& \ddots & &\vdots \\
		& &\ddots& &\ddots & \ddots&0\\
		& &  & \theta&\cdots&0 & 1\\
		& &  & & \theta&\cdots& 0\\
		& &  & & &\ddots& \vdots\\
		& &  & & & & \theta\\
		\end{pmatrix}+\begin{pmatrix}
		0&\cdots&\cdots & \cdots &\cdots & \cdots & 0\\
		\vdots& & & & & &\vdots\\
		0& & & & & & 0\\
		1&0 &\cdots &\cdots &\cdots &\cdots &0\\
		&\ddots& \ddots&  & & &\vdots\\
		& &1& \ddots & & & \vdots\\
		k_1&\cdots&\cdots &k_r&0&\cdots& 0
		\end{pmatrix}\tau.
		\]
In particular, the last $r$-coordinates of $\rho^{[k]}$ are tractable. By convention, we also let $\rho^{[0]}=\phi$.

For any $f\in \mathbb{F}_q(t)$, we define 
\[
\dmat_k[f]:=\begin{pmatrix}
		f& &\partial_{t}(f)& & \partial_{t}^2(f)  & & \partial^k_{t}(f) \\
		&\ddots  & &\ddots& &  \ddots &\vdots\\
		& & \ddots&  &\ddots & &\partial_{t}^2(f) \\
		& &  &\ddots & &\ddots &\\
		& &  & &\ddots & & \partial_{t}(f)\\
		& &  & & & \ddots& \\
		& &  & & & & f\\
		\end{pmatrix}|_{t=\theta}\in \Mat_{rk+1}(K).
        \]
Indeed, in the first row of $\dmat_k[f]$, $\partial_{t}^i(f)$ appears in the $(ir+1)$-st entry for $0\leq i \leq k$ and the remaining entries are zero. Moreover, when $f\in \bA$, $\dmat_k[f]$ is the constant coefficient of $\rho^{[k]}_{f}$.

By \eqref{E:RAT1} and \eqref{E:AGF4}, we note that $\omega_C^k\Upsilon_{\phi}\in \GL_r(\mathbb{T})$ is a rigid analytic trivialization of $M_{k}$ and hence $M_k$ is rigid analytically trivial. Thus, by \cite[Thm. 4]{And86}, $\rho^{[k]}$ is uniformizable and its period lattice $\Lambda_{\rho^{[k]}}$ is of rank $r$ over $\bA$.

Next, we assume that $L$ is a perfect field and analyze the $t$-comotive $N_k$ of $\rho^{[k]}$ which is the $L[t,\sigma]$-module given by
$
N_k:=N_{\phi}\otimes_{L[t]}N_{C^{\otimes k}}.
$
Then one can choose an $L[t]$-basis $\mathfrak{B}_k:=[1\otimes \bn_k,\dots,\sigma^{r-1}\otimes \bn_k]^{\tr}\in \Mat_{r\times 1}(N_k)$ for $N_k$ so that 
$
\sigma \cdot \mathfrak{B}_k=(t-\theta)^k\Phi_{\phi} \mathfrak{B}_k.
$
Moreover, if we set 
$
\Psi_k=\Omega(t)^k\Psi_{\phi}=(\omega_C^{-k})^{(1)}\Psi_{\phi}\in \GL_r(\mathbb{T}),
$ then we see that 
$
\Psi_k^{(-1)}=(t-\theta)^k\Phi_{\phi}\Psi_k
$
and hence $N_{k}$ is rigid analytically trivial.

\subsection{The period lattice of $\rho^{[k]}$ and the map $\varphi_{\rho^{[k]}}$} In what follows, we will use the methods in \cite[\S5]{Mau22} developed by Maurischat to describe the entries of the generators of the period lattice $\Lambda_{\rho^{[k]}}$ of $\rho^{[k]}$. Let $M^1_k$ be the $L[t]$-submodule of $M_k$ generated by $\tau\cdot \mathfrak{m}_1\otimes \bm,\dots, \tau\cdot \mathfrak{m}_r\otimes \bm$. We note, by \eqref{E:extmapk}, that  
\[\mathfrak{F}^1_{k}:=[(t-\theta)^k\mathfrak{m}_1\otimes \bm,(t-\theta)^{k}\mathfrak{m}_2\otimes \bm,\dots, (t-\theta)^{k}\mathfrak{m}_r\otimes \bm]^{\tr}\in \Mat_{r\times 1}(M^1_k)\]
forms an $L[t]$-module basis for $M^1_k$.

For $1\leq \nu \leq r-1$, we define \textit{the quasi-periodic function $F_{\tau^{\nu}}:\mathbb{C}_{\infty}\to \mathbb{C}_{\infty}$ of $\phi$ with respect to the biderivation that maps $\theta$ to $\tau^{\nu}$} to be the unique function satisfying
\begin{equation}\label{D:quasi periodic function}
F_{\tau^{\nu}}(\theta z)-\theta F_{\tau^{\nu}}(z)=\Exp_{\phi}(z)^{q^{\nu}}, \ \ \ \ z\in \mathbb{C}_{\infty}.
\end{equation}
We refer the reader to \cite[\S4.1]{NP21} for more details on quasi-periodic functions and biderivations. 

We let 
\[
A':=\begin{pmatrix}
    \frac{1}{(t-\theta)^{k}} & & & &  \\
    & \ddots & & &  \\
    &  & \ddots & &  \\
    & & & &  \frac{1}{(t-\theta)^{k}}  \\
   \frac{1}{(t-\theta)^{k-1}} & & & &  \\
    & \ddots & & &  \\
    &  & \ddots & &  \\
    & & & &  \frac{1}{(t-\theta)^{k-1}}  \\
    \vdots & & &  & \\
    \vdots &  & & &    \\
    \frac{1}{(t-\theta)} & & & &  \\
    & \ddots & & &  \\
    &  & \ddots & &  \\
    & & & &  \frac{1}{(t-\theta)}  \\
    1& &  & &      \end{pmatrix}\in \Mat_{(rk+1)\times r}(K[t]).
\]
Recall the $L[\tau]$-basis $\bb_k$ of $M_k$ defined in \S\ref{S:drinfeldtens}. Note that $\bb_k=A'\mathfrak{F}^1_k$. For any $\alpha_1,\dots,\alpha_r\in L[t]$, let $\diag(\alpha_1,\dots,\alpha_r)\in \Mat_r(\mathbb{T})$ be the diagonal matrix whose $i$-th diagonal entry is $\alpha_i$. 
Then, by \cite[Cor. 5.5, Rem. 5.7]{Mau22}, for $1\leq j \leq r$, a generator $\varpi_{k,j}\in \mathbb{C}_{\infty}^{rk+1}$ of $\Lambda_{\rho^{[k]}}$ may be chosen as the $j$-th column of $\Res_{t=\theta}(A'\diag((t-\theta)^{k},(t-\theta)^{k},\dots,(t-\theta)^{k})\omega_C^k\Upsilon_{\phi})$. Recall the Anderson generating function $f_{\lambda_j}(t)$ of $\phi$ with respect to $\lambda_j$. Observe that the $j$-th column of $ A'\diag((t-\theta)^{k},(t-\theta)^{k},\dots,(t-\theta)^{k})\omega_C^k\Upsilon_{\phi} $ is given by 
\begin{multline*}
\mathcal{F}_{k,j}:=\omega_C^k(f_{\lambda_j}(t),
\dots,
f_{\lambda_j}(t)^{(r-1)},
(t-\theta)f_{\lambda_j}(t),\dots
,(t-\theta)f_{\lambda_j}(t)^{(r-1)},\dots,\\(t-\theta)^{k-1}f_{\lambda_j}(t),\dots,(t-\theta)^{k-1}f_{\lambda_j}(t)^{(r-1)}
,
(t-\theta)^{k}f_{\lambda_j}(t)
)\in \mathbb{T}^{rk+1}.
\end{multline*}
Thus, using  the fact that $\Res_{t=\theta}\omega_C=-\tilde{\pi}$, we can choose $\varpi_{k,j}$  so that
\begin{equation}\label{E:periodsofext}
\varpi_{k,j}=\Res_{t=\theta}\mathcal{F}_{k,j}=(*,\dots,*,-\tilde{\pi}^{k}F_{\tau}(\lambda_j),\dots,-\tilde{\pi}^{k}F_{\tau^{r-1}}(\lambda_j),\tilde{\pi}^k\lambda_j)^{\tr}\in \mathbb{C}_{\infty}^{rk+1}.
\end{equation}
Thus, we have shown the following.
\begin{proposition}\label{P:periodbasis}
    For any $k\geq 1$, the set $\{\varpi_{k,1},\dots,\varpi_{k,r}\}$ forms an $\bA$-module basis for $\Lambda_{\rho^{[k]}}$.
\end{proposition}

\begin{remark} To apply \cite[Cor. 5.5, Rem. 5.7]{Mau22}, Maurischat assumes Anderson modules to be defined over a finite extension of $K$. However, we emphasize that the exact same idea of the proof of these results work for Anderson modules defined over any subfield of $\mathbb{C}_{\infty}$ containing $K$. We would like to thank A. Maurischat for clarifying this with the authors. 
\end{remark}

We continue with an analysis on the entries of Anderson generating functions of $\rho^{[k]}$ with respect to the generators of the period lattice $\Lambda_{\rho^{[k]}}$. For each $1\leq j \leq r$, we consider the Anderson generating function $\mathcal{Y}_{\varpi_{k,j}}(t)$ corresponding to $-\varpi_{k,j}$ which is given by 
\[
\mathcal{Y}_{j}(t):=\mathcal{Y}_{\varpi_{k,j}}(t)=\sum_{i=0}^{\infty}\Exp_{\rho^{[k]}}\left((\dmat_k[t])^{-i-1}(-1)\varpi_{k,j}\right)t^i=\begin{pmatrix}
    \mathcal{Y}_{1,j}\\
    \vdots\\
   \mathcal{Y}_{rk+1,j}
\end{pmatrix}\in \mathbb{T}^d.
\]
By Lemma \ref{L:AGF}(iii), we obtain 
\begin{equation}\label{E:AGF1}
\rho^{[k]}_{t}(\mathcal{Y}_{\varpi_{k,j}}(t))=t\mathcal{Y}_{\varpi_{k,j}}(t).
\end{equation}

\begin{lemma}\label{L:coef} For $0\leq u \leq k-1$ and $1\leq \mu \leq r$, we have 
\[
\mathcal{Y}_{ru+\mu,j}=(t-\theta)^{u}\frac{\mathcal{Y}_{1,j}^{(\mu-1)}}{\mathbb{S}_{\mu-2}^k}
\]
and 
\[
\mathcal{Y}_{rk+1,j}=(t-\theta)^{k}\mathcal{Y}_{1,j}.
\]
Furthermore, we have $\mathcal{Y}_{1,j}=\omega_C^kf_{\lambda_j}(t)$. 
\end{lemma}
\begin{proof} The first assertion follows from a direct computation using \eqref{E:AGF1}. The second assertion is a consequence of the first assertion and comparing the last entry of both sides of \eqref{E:AGF1}.  Finally, the last assertion follows from the first assertion, \cite[Prop. 3.8]{Mau22} and the fact that $\Res_{t=\theta}\mathcal{F}_{k,j}=\varpi_{k,j}$.  
\end{proof} 
\begin{remark}\label{R:rat}Note that by Lemma \ref{L:coef}, we see that $\Psi_{k}$ is exactly the same as the matrix $\Psi$ described in \cite[\S2.4.1]{GG26} for $\rho^{[k]}$. This observation is crucial to state our formulas for $\rho^{[k]}$ purely in terms of the rigid analytic trivialization of the Drinfeld module $\phi$ in Theorem \ref{T:mellin}.
\end{remark}

We finish this subsection by describing the map $\varphi_{\rho^{[k]}}$. Define 
\begin{multline}\label{D:M phi2}
\mathbb{M}_{\rho^{[k]}}=\{\left(\sum_{n=0}^{\infty}a_{1,n}\tau^n,\dots,\sum_{n=0}^{\infty}a_{rk+1,n}\tau^n\right)\in \Mat_{1\times (rk+1)}(\mathbb{C}_{\infty}[[\tau]])\ \ | \ \  \\ \text{ for each } 1\leq \ell \leq rk+1\ \ \text{ and }1\leq \nu \leq r, \ \  
|a_{\ell,n}|||\mathcal{Y}_{\nu,\ell}^{(n)}||\to 0 
\text{ as } n\to \infty \}.
\end{multline}
We also let $\mathfrak{S}_{-1}=\Id_r$ and for $n\geq 0$, set $\mathfrak{S}_{n}=\mathbb{S}_n (\Phi_{\phi}^{\tr})^{(n)}\cdots (\Phi_{\phi}^{\tr})^{(1)}\Phi_{\phi}^{\tr} \in \GL_r(\mathbb{C}_{\infty}(t))\cap\GL_r(\mathbb{T}) $.  The map
$\varphi_{\rho^{[k]}}:\mathbb{M}_{\rho^{[k]}}\to \Mat_{1\times r}(\mathbb{T})$ is given by  
\begin{multline*}
\varphi_{\rho^{[k]}}\left(\left(\sum_{n= 0}^{\infty}a_{1,n}\tau^n,\dots,\sum_{n= 0}^{\infty}a_{rk+1,n}\tau^n\right)\right)=\sum_{\substack{0\leq u \leq k-1\\1\leq j \leq r}}\sum_{n= 0}^{\infty}a_{ru+j,n}\mathfrak{e}_j^{\tr}\mathbb{S}_{n-1}^k(t-\theta^{q^n})^u\mathfrak{S}_{n-1}\\
+\sum_{n=0}^{\infty}a_{rk+1,n}\mathfrak{e}_1^{\tr}\mathbb{S}_{n}^k\mathfrak{S}_{n-1}.
\end{multline*}

\section{Dual Goss \texorpdfstring{$L$}{L}-functions  and proof of Theorem \ref{T:1}}\label{S:3}
In this section, our main goal is to prove Theorem \ref{T:1}. In the next subsection, we provide necessary details for the construction of dual Goss $L$-functions as well as Taelman $L$-values and we refer the reader to  \cite{Fang}, \cite[\S5]{Gek91} and \cite{Gos92} for further details. 
\subsection{Dual Goss \texorpdfstring{$L$}{L}-functions of Drinfeld modules and Taelman $L$-values} \label{S:lfunctionsTaelman}
Let $A_{+}$ be the set of monic polynomials in $A$. Consider an irreducible polynomial $\mathfrak{p}\in A_{+}$ and set $L:=A/\mathfrak{p} A$.  
Let us set $\overline{\phi}$ be the Anderson module where $\overline{\phi}_a$ is the reduction of $\phi_{a}$ modulo $\mathfrak{p}$ for each $a\in\bA$. 
Let $M$ be the $t$-motive corresponding to $\phi$ which has rank $r$ over $K[t]$. We say that $M$ has \textit{good reduction at $\mathfrak{p}$} if the  $t$-motive $M_{\overline{\phi}}$ over $L$ corresponding to the abelian Anderson module $\overline{\phi}$  has rank $r$ over $L[t]$. Here, by $M_{\overline{\phi}}$, we mean the $L[t,\tau]$-module $M_{\overline{\phi}}=\Mat_{1\times d}(L[\tau])$ with an $L[t]$-module structure given by 
\[
\alpha t^i\cdot m:=\alpha  m  \overline{\phi}_{t^i}, \ \ m\in M_{\overline{\phi}}, \ \ \alpha\in L.
\]
    We further say that $M$ has \textit{everywhere good reduction} if $M$ has good reduction at any irreducible $\mathfrak{p}\in A_{+}$.
		
		Let $K^{\text{sep}}$ be the separable closure of $K$ in $\CC_{\infty}$ and define $M_{K^{\text{sep}}}:=M\otimes_{K} K^{\text{sep}}$. Let $v$ be a monic irreducible element in $\textbf{A}$. We define  $\textbf{A}_v:=\varprojlim_{i}\textbf{A}/v^i\textbf{A}$ and consider $\textbf{K}_v:=\textbf{A}_v\otimes_{\textbf{A}}\mathbb{F}_q(t)$. We set
		\[
		\mathcal{H}_{v}(M):=\varprojlim_{i}(M_{K^{\text{sep}}}/v^iM_{K^{\text{sep}}})^{\tau}
		\]
		where, for each $i$,  $(M_{K^{\text{sep}}}/v^iM_{K^{\text{sep}}})^{\tau}$ denotes the set of elements in  $M_{K^{\text{sep}}}/v^iM_{K^{\text{sep}}}$ fixed by the action of $\tau$. Note that $\mathcal{H}_{v}(M)$ is a free $\bA_v$-module of rank $r$ with a continuous action of $\Gal(K^{\text{sep}}/K)$. Moreover, by \cite[Thm. 1.1, Thm. 7.3]{G01}, for any irreducible element $v\in \bA$ so that $v(\theta)$ is not equal to $\mathfrak{p}$,  $I_{\mathfrak{p}}$ acts trivially on $\mathcal{H}_{v}(M)$ if and only if $M$ has good reduction at $\mathfrak{p}$ (see \cite[\S5.6]{Gos96} for more details). 
		
		Let $\mathfrak{p}'$ be a prime of $K^{\text{sep}}$ lying above $\mathfrak{p}$ and  $I_{\mathfrak{p}'}\subset \Gal(K^{\text{sep}}/K)$ be the inertia group at $\mathfrak{p}'$. We set $\mathcal{H}_{v}(M)^{I_{\mathfrak{p}'}}$ to be the $\bA_v$-module consisting of elements of $\mathcal{H}_{v}(M)$ invariant under the action of $I_{\mathfrak{p}'}$. Let $\overline{L}$ be a fixed algebraic closure of $L$ and let $\Frob_{\mathfrak{p}}\in \Gal(\overline{L}/L)$ be the arithmetic Frobenius at $\mathfrak{p}$. We define \textit{the local factor of $M$ at $\mathfrak{p}$} by
		\[
		P_{\mathfrak{p}}^{M}(X):=\det_{\bA_v}(1-X\Frob_{\mathfrak{p}}\ \ | \ \ \mathcal{H}_{v}(M)^{I_{\mathfrak{p}'}})\in \bA_v[X].
		\]
		Note that the definition of $P_{\mathfrak{p}}^{M}(X)$ is indeed independent of the choice of ${\mathfrak{p}'}$ \cite[\S8.6]{Gos96}. 
       By \cite[Thm. 1.1, Thm. 7.3]{G01}, we further know that for such $v$ and $\mathfrak{p}$, the polynomial $P_{\mathfrak{p}}^{M}(X)$ has coefficients in $\mathbb{F}_q(t)$ and is independent of the choice of $v$. 
        
		There exists a finite set $S$ of primes in $A$ where $M$ does not have good reduction. For any non-negative integer $k$, consider 
		\[
L_S(M^{\vee},s):=\prod_{\mathfrak{p}\not\in S}P^{M}(\mathfrak{p}_{|\theta=t}^{-s})^{-1}.
		\]
By \cite[Cor. 8.4(ii), (iii)]{G01}, we have $L_S(M^{\vee},s)\in \mathbb{F}_q(t)^{\times}$. We finally define \textit{the dual Goss $L$-function $L(M^{\vee},k)$ of $M$} by
    \[
    L(M^{\vee},s):=L_S(M^{\vee},s)\prod_{\mathfrak{p}\in S}P^{M}(\mathfrak{p}_{|\theta=t}^{-s})^{-1}\in \mathbb{F}_q((1/t))^{\times}.
    \]
By convention, when $M_{\phi}$ is the abelian $t$-motive corresponding to a Drinfeld module $\phi$, we further write
\[
L(\phi^{\vee},s):=L(M_{\phi}^{\vee},s)|_{t=\theta}\in K_{\infty}.
\]
As an example, the value of the dual Goss $L$-function $L(C^{\vee},s)=L_{\emptyset}(M_{C}^{\vee},s)$ at $s$ corresponding to the $t$-motive $M_{C}$ is given by 
		\[
		L(C^{\vee},s)=\sum_{a\in A_{+}}\frac{1}{a^{s+1}}\in K_{\infty}^{\times},
		\]
which are so called \textit{the Carlitz zeta values} and they can be seen as a function field analogue of Riemann zeta values. 
Furthermore, when $\phi$ is a Drinfeld module defined over $A$, by \cite[Cor. 3.6(a)]{CGP18}, we know that $L(\phi^{\vee},s)$ converges in $K_{\infty}$ for any non-negative integer $s$.

		In what follows, we introduce Taelman $L$-values. Let $\phi$ be an abelian Anderson module of dimension $d$ defined over $A$. For any finite $\bA$-module $\mathrm{N}$ isomorphic to the direct sum $\oplus_{i=1}^s\bA/g_i\bA$ for some monic $g_i\in \bA$, we set $|\mathrm{N}|_{A}:=\prod_{i=1}^sg_i|_{t=\theta}\in A\setminus \{0\}$. Then \textit{the Taelman $L$-value $L(\phi/A)$} is defined by the infinite product
		\[
		L(\phi/A):=\prod_{\substack{\mathfrak{p}\in A_{+}\\ \mathfrak{p} \text{ irreducible}}}\frac{|\Lie(\overline{\phi})(A/\mathfrak{p} A)|_{A}}{|\overline{\phi}(A/\mathfrak{p} A)|_{A}}\in 1+\frac{1}{\theta}\mathbb{F}_q\Big[\Big[\frac{1}{\theta}\Big]\Big].
		\]
\begin{proposition}\label{P:lvalueTaelman} Let $\phi$ be a Drinfeld module of rank $r$ defined over $A$ and $k$ a non-negative integer. Let $\rho^{[k]}$ be the abelian Anderson module defined over $A$ given as in \S\ref{S:drinfeldtens}. Then we have 
                \[
        L(\phi^{\vee},k)=\mathfrak{c}_{k}L(\rho^{[k]}/A)
                \]
                    for some $\mathfrak{c}_k\in K^{\times}$. Moreover, if $\phi$ has everywhere good reduction, then $\mathfrak{c}_{k}=1$.
                \end{proposition}
                \begin{proof}
                  Recall from \S\ref{S:drinfeldtens} that $M_k=M_{\phi}\otimes_{K[t]}M_{C^{\otimes k}}$ and $\rho^{[k]}$ is the abelian Anderson module corresponding to $M_k$. Proceeding as in \cite[\S3.3.3]{GN24}, we see that $L_S(M_k^{\vee},0)$ differs from $L(\rho^{[k]}/A)|_{\theta=t}$ by an element in $\mathbb{F}_q(t)^{\times}$. On the other hand, since $L(M_k^{\vee},0)$ also differs from $L_S(M_k^{\vee},0)$ by an element in $\mathbb{F}_q(t)^{\times}$,  we can write 
                  \begin{equation}\label{E:Lanalysis1}
                  L(M_k^{\vee},0)=\alpha_1 L(\rho^{[k]}/A)|_{\theta=t}
                  \end{equation}
                  for some $\alpha_1 \in \mathbb{F}_q(t)^{\times}$. Moreover, by \cite[Thm. 2.2]{Ye26}, we see that 
                  \begin{equation}\label{E:Lanalysis2}
                  L_S(M_k^{\vee},0)=L_S(\phi^{\vee},k).
                  \end{equation}
                  Combining \eqref{E:Lanalysis1} and \eqref{E:Lanalysis2} yields 
                  \[
                  L(\phi^{\vee},k)=\alpha_2 L(\rho^{[k]}/A)                 \]
                  for some $\alpha_2\in K^{\times}$. Hence, we obtain the first assertion. Since, when $\phi$ has everywhere good reduction, we have $S=\emptyset$, the second assertion follows from the fact that $L_S(M_k^{\vee},0)=L(M_k^{\vee},0)=L(\phi^{\vee},k)=L(\rho^{[k]}/A)|_{\theta=t}$.
                  \end{proof}
                      
\subsection{Proof of Theorem \ref{T:1}} \label{S:analysisunit}
Using Lemma \ref{L:coef}, Remark \ref{R:rat} and \cite[Thm. 1.6]{GG26}, we first record the Mellin transform formulas for the logarithm function of $\rho^{[k]}$. For each $1\leq \ell \leq rk+1$, recall the element $\mathfrak{N}_{\ell}\in \mathbb{M}_{\rho^{[k]}}$, the matrix $\mathcal{F}_k\in \Mat_{(rk+1)\times r}(\mathbb{T})$ as well as the integer $\mathfrak{m}_{k}$ from \S\ref{S:mainresults}. Recall also that for $k\geq 0$, we define
\[
\Psi_k=\frac{1}{(\omega_C^{(1)})^k}\Psi_{\phi}\GL_r(\mathbb{T})
\]
where $\Psi_{\phi}\in \GL_r(\mathbb{T})$ is as defined in \eqref{E:drinfeldpsi}.
\begin{theorem}\label{T:mellin} Let $\mathfrak{p}_{\ell}:\mathbb{C}_{\infty}^{rk+1}\to \mathbb{C}_{\infty}$ be the projection onto the $\ell$-th coordinate.  Then, for any $\zz\in \mathbb{C}_{\infty}^{rk+1}$ lying in the domain of convergence of $\Log_{\rho^{[k]}}\theta^{-\mathfrak{m}_k}$, we have
\begin{equation}\label{E:formula}
\mathfrak{p}_{\ell}(\Log_{\rho^{[k]}}(\theta^{-\mathfrak{m}_k}\zz))=
\mathcal{M}_{\zz}\left(\mathfrak{N}_{\ell}\cdot \mathcal{F}_k((\Psi_k)^{\tr})^{(-1)}\right).
\end{equation}
Moreover, let $r(k-1)+2\leq \mathfrak{t}\leq rk+1$ and for each $1\leq \nu \leq r$, consider $\varpi_{k,\nu}=(\varpi_{k,\nu,1},\dots,\varpi_{k,\nu, rk+1})^{\tr}\in \mathbb{C}_{\infty}^{rk+1}$  given as in \eqref{E:periodsofext}. If $\mathfrak{m}_k=0$ and for $i\geq 1$, each $\dmat_{k}[t]^{-i}\varpi_{k,\nu}$ lies in the domain of convergence of $\Log_{\rho^{[k]}}$, then we have 
\[
\mathfrak{p}_{\mathfrak{t}}(\Log_{\rho^{[k]}}(\zz))=\mathcal{M}_{\zz}\left(\frac{1}{\theta-t}[\varpi_{k,1, \mathfrak{t}},\dots, \varpi_{k,r, \mathfrak{t}}](\Psi_{k}^{\tr})^{(-1)}\right).
\]

\end{theorem}

Let $V$ be an $\mathbb{F}_q((1/t))$-vector space of dimension $r$. By \textit{an $\bA$-lattice of rank $r$}, we mean a free and finitely generated $\bA$-module $V'$ of rank $r$ in $V$ containing an $\mathbb{F}_q((1/t))$-basis of $V$. Let $V_1$ and $V_2$ be two $\bA$-lattices in $V$.  By \cite[Lem. 1]{AR16}, we know that the rank of $V_1$ and $V_2$ as $\bA$-modules are the same as the dimension of $V$ as an $\mathbb{F}_q((1/t))$-vector space.  Let us choose such $\mathbb{F}_q((1/t))$-vector space bases $\{\bv_{11},\dots, \bv_{1r}\}$ and $\{\bv_{21},\dots,\bv_{2r}\}$ in $V_1$ and $V_2$ for $V$ respectively and consider the $\mathbb{F}_q((1/t))$-linear map
\[
f:V=\bigoplus_{i=1}^{r} \mathbb{F}_q((1/t))\cdot \bv_{1i}\to V=\bigoplus_{i=1}^{r} \mathbb{F}_q((1/t))\cdot \bv_{2i}
\]
defined by $f(\bv_{1i})=\bv_{2i}$ for each $1\leq i \leq r$. For each $1\leq i \leq r$, there exist $\bb_{i1},\dots,\bb_{ir}\in \mathbb{F}_q((1/t))$ such that
\begin{equation}\label{E:reg1}
\bv_{2i}=\bb_{i1}\cdot\bv_{11}+\cdots+\bb_{ir}\cdot \bv_{1r}.
\end{equation}

Set $\mathfrak{D}:=(\bb_{ij})_{i,j}\in \GL_r(\mathbb{F}_q((1/t)))$. We further call any $g=\sum_{i\geq i_0}c_it^{-i}\in \mathbb{F}_q((1/t))^{\times}$ \textit{monic} if $c_{i_0}=1$ and then define
\[
[V_1:V_2]:=\mathfrak{c}\det(\mathfrak{D}).
\]
where $\mathfrak{c}\in \mathbb{F}_q^{\times}$ so that $\mathfrak{c}\det(\mathfrak{D})$ is monic.

Let $\phi$ be a Drinfeld module of rank $r$ defined over $A$ and let $\rho^{[k]}$ be the abelian Anderson module defined over $A$ given as in \S\ref{S:drinfeldtens}. Next, we define the $\mathbb{F}_q((1/t))$-action on $\Lie(\rho^{[k]})(K_{\infty})$ by
\[
\mathfrak{y}\cdot x:=\sum_{i\geq i_0}c_i\dmat_k[t]^{-i}x\in \Lie(\rho^{[k]})(K_{\infty}), \ \ \  \mathfrak{y}=\sum_{i\geq i_0}c_it^{-i}\in \mathbb{F}_q((1/t)), \ \ \ \ x\in \Lie(\rho^{[k]})(K_{\infty}).
\]
Let $W$ be the $\mathbb{F}_q((1/t))$-vector space in $\Lie(\rho^{[k]})(K_{\infty})$ constructed in \cite[Prop. 4.3]{ANDTR20} which is isomorphic to the quotient 
\[
W\cong \Lie(\rho^{[k]})(K_\infty)/ (\dmat_k -\theta \Id_{rk+1})\Lie(\rho^{[k]})(K_\infty)
\]
via the natural projection map as $\mathbb{F}_q((1/t))$-vector spaces (see Proposition \ref{P:intersection} for the construction of such $W$ in the case of a particular abelian Anderson module). Since the last $r$ coordinates of $\rho^{[k]}$ are tractable, we see that $W$ is an $r$-dimensional $\mathbb{F}_q((1/t))$-vector space in $\Lie(\rho^{[k]})(K_{\infty})$. \textit{The unit module of $\rho^{[k]}$} is given by 
\[
U(\rho^{[k]}/A):=\left\{ x\in \Lie(\rho^{[k]})(K_{\infty}) \ \ | \ \ \Exp_{\rho^{[k]}}(x)\in \rho^{[k]}(A) \right\}.
\]
By \cite[Thm. 4.4]{ANDTR20}, we know that the $\bA$-modules $V_1:=\Lie(\rho^{[k]})(A)\cap W$ and $V_2:=U(\rho^{[k]}/A)\cap W$ are indeed $\bA$-lattices of rank $r$ in $W$.

Our goal now is to analyze $[V_1:V_2]$. As above, let $\{\bv_{1,1},\dots, \bv_{1,r}\}$ and $\{\bv_{2,1},\dots,\bv_{2,r}\}$ be $\mathbb{F}_q((1/t))$-vector space bases in $V_1$ and $V_2$ for $W$ respectively. For each $1\leq i \leq r$ and $j=1,2$, let us write $\bv_{j,i}=(\bv_{j,i,1},\dots,\bv_{j,i,rk+1})^{\tr}$. Since the last $r$ coordinates of $\rho^{[k]}$ are tractable, \eqref{E:reg1} implies 
\begin{equation}\label{E:tractablereg}
\bv_{2,i,r(k-1)+\mu+1}=\bb_{i1}\cdot\bv_{1,1,r(k-1)+\mu+1}+\cdots+\bb_{ir}\cdot \bv_{1,r,r(k-1)+\mu+1}, \ \ \ \ \ 1\leq \mu \leq r.
\end{equation}
Then, setting $\bV_2:=(\bv_{2,i,(r(k-1)+\mu+1)}|_{\theta=t})_{i,\mu}\in \Mat_r(\mathbb{F}_q((1/t)))$ and $\bV_1:=(\bv_{1,i,r(k-1)+\mu+1}|_{\theta=t})_{i,\mu}\in \Mat_r(\bA)$, we obtain 
\[
\bV_2=\bV_1\mathfrak{D}.
\]
Thus, by \cite[Thm. 4.4]{ANDTR20}, we see that $\bV_1\in \GL_r(\mathbb{F}_q((1/t))$, $\bV_2\in \GL_r(\mathbb{F}_q(t))$ and 
\begin{equation}\label{E:lvalueformula}
L(\rho^{[k]}/A)=\alpha\det(\mathfrak{D})|_{t=\theta}=\alpha\det(\bV_1)|_{t=\theta}^{-1}\det(\bV_2)_{t=\theta}=\alpha'\det(\bV_2)|_{t=\theta}
\end{equation}
for some $\alpha, \alpha'\in K^{\times}$. For each $1\leq i \leq r$, we let $\widetilde{\bv_{2,i}}\in K_{\infty}^{rk+1}$ be the vector obtained by evaluating each entry of $\bv_{2,i}$ at $t=\theta$. Then, we see that, up to a constant in $K^{\times}$, $L(\rho^{[k]}/A)$ is the determinant of the $r\times r$-matrix whose $i$-th row is the last $r$ entries of $\widetilde{\bv_{2,i}}$.

\begin{proof}[{Proof of Theorem \ref{T:1}}] Observe that \eqref{E:detlvalue2} is already shown in \cite[Thm. 4.2]{GG26}. We now let $k\geq 1$ and show that \eqref{E:detlvalue1} holds. Since, by \cite[Lem. 2.5.4]{HJ20}, $\Exp_{\rho^{[k]}}$ is an isometric map with its inverse $\Log_{\rho^{[k]}}$ near $0$, for each $1\leq i \leq r$, there exists a non-negative integer $N_k$ such that $\Exp_{\rho^{[k]}}(\theta^{-N_k}\widetilde{\bv_{2,i}})\in \Lie(\rho^{[k]})(K_{\infty})$ lies in the domain of convergence of $\Log_{\rho^{[k]}}$. We set 
\[
\zz_i:=\theta^{\mathfrak{m}_k}\Exp_{\rho^{[k]}}\left(\frac{\widetilde{\bv_{2,i}}}{\theta^{N_k}}\right)\in \mathbb{C}_{\infty}^{rk+1}.
\]
Thus, by Theorem \ref{T:mellin}, for $1\leq \mu \leq r$, the $r(k-1)+\mu+1$-th entry of $\widetilde{\bv_{2,i}}$ is given by $\theta^{N_k}\mathcal{M}_{\mathfrak{N}_{r(k-1)+\mu+1},\zz_{i}}\left(\mathcal{F}_k((\Psi_k)^{\tr})^{(-1)}\right)$. 
  Then \eqref{E:detlvalue1} follows from Proposition \ref{P:lvalueTaelman} and \eqref{E:lvalueformula}.

What remains to show that the elements $\zz_1,\dots,\zz_r\in \mathbb{C}_{\infty}^{rk+1}$ ($\zz\in \mathbb{C}_{\infty}$ for $k=0$ respectively) lie in $\overline{K}^{rk+1}$ (in $\overline{K}$ respectively). Let us first assume that $k=0$. In this case, by \cite[Thm. 4.2(i)]{GG26}, one can choose $\zz=\theta^{\mathfrak{m}_{\phi}}\Exp_{\phi}(\mathfrak{f}/\theta^{N})$ for some non-negative integer $N$ and some $\mathfrak{f}\in K_{\infty}$ such that $\Exp_{\phi}(\mathfrak{f})\in A$. Then, we have
\[
\phi_{t^N}\left(\Exp_{\phi}(\mathfrak{f}/\theta^{N})\right)=\Exp_{\phi}(\mathfrak{f})\in A.
\]
    In other words, $\Exp_{\phi}(\mathfrak{f}/\theta^{N})$ is a root of a polynomial of degree $Nq^{r}$ with coefficients in $\overline{K}$ so that its constant term is an element in $A$. Hence $\Exp_{\phi}(\mathfrak{f}/\theta^{N})\in \overline{K}$, implying that $\zz\in \overline{K}$. 
    
    Now let $k\geq 1$. Recall that for each $1\leq i \leq r$, we have $\Exp_{\rho^{[k]}}(\widetilde{\bv_{i}})\in A^{rk+1}$. Observe that 
\begin{equation}\label{E:claimapply}
\rho^{[k]}_{t^{N_k}}\left(\Exp_{\rho^{[k]}}\left(\dmat_k[t^{N_k}]^{-1}\widetilde{\bv_{i}}\right)\right)=\Exp_{\rho^{[k]}}\left(\widetilde{\bv_{i}}\right)\in A^{rk+1}.
\end{equation}
Let $n\geq 1$. We claim that  if $X=(x_1,\dots,x_{rk+1})^{\tr}\in \mathbb{C}_{\infty}^{rk+1}$ satisfies \begin{equation}\label{E:soe11}
\rho^{[k]}_{t^n}(X)=(a_1,\dots,a_{rk+1})^{\tr}\in \overline{K}^{rk+1},
\end{equation}
then $X\in \overline{K}^{rk+1}$. We proceed by induction on $n$. When $n=1$, by the construction of $\rho^{[k]}_t$, we see that comparing the first $r(k-1)+1$-st entry of both sides of \eqref{E:soe11} implies
\begin{equation}\label{E:soe12}
\theta x_{ru+j}+x_{r(u+1)+j}=a_{ru+j}, \ \ 0\leq u \leq k-1, \ \ 1\leq j \leq r.
\end{equation}
On the other hand, comparing the next $(r-1)$-st entry of both sides of \eqref{E:soe11} gives
\begin{equation}\label{E:soe13}
    x_\nu^q+\theta x_{r(k-1)+1+\nu}=a_{r(k-1)+1+\nu}, \ \ 1\leq \nu \leq r-1.
\end{equation}
Finally, the comparison of the last entries of both sides of \eqref{E:soe11} yields
\begin{equation}\label{E:soe14}
k_1x_1^q+\cdots+k_rx_r^q+\theta x_{rk+1}=a_{rk+1}.
\end{equation}
Using \eqref{E:soe12}, one can write each entry of $X$ in terms of a $\overline{K}$-linear combination of $1,x_1,\dots,x_r$. Combining this with \eqref{E:soe13}, one can then write each entry of $X$ in terms of a polynomial of $x_1$ with coefficients in $\overline{K}$. Thus, \eqref{E:soe14} yields that $x_1$ is a solution of a polynomial of order $q^r$ with coefficients in $\overline{K}$. Therefore, $x_1\in \overline{K}$, implying, by the above discussion, that $X\in \overline{K}^{rk+1}$. Now assume that the claim holds for $1\leq i < n$. We write 
\[
\rho^{[k]}_{t^n}(X)=\rho^{[k]}_t(\rho^{[k]}_{t^{n-1}}(X))=(a_1,\dots,a_{rk+1})^{\tr}
\]
and set $Y:=\rho^{[k]}_{t^{n-1}}(X)\in \mathbb{C}_{\infty}^{rk+1}$. Since, we have $\rho^{[k]}_{t}(Y)\in \overline{K}^{rk+1}$, by the induction hypothesis, we see that $Y\in \overline{K}^{rk+1}$. Thus, since now we have $Y=\rho^{[k]}_{t^{n-1}}(X)\in \overline{K}^{rk+1}$, again by the induction hypothesis, we obtain $X\in \overline{K}^{rk+1}$ as desired. Thus, applying our claim to the equality \eqref{E:claimapply}, we conclude that $\Exp_{\rho^{[k]}}\left(\dmat_k[t^{N_k}]^{-1}\widetilde{\bv_{i}}\right)\in \overline{K}^{rk+1}$. Hence $\zz_i\in \overline{K}^{rk+1}$, finishing the proof of the theorem.
\end{proof}
\section{Logarithm of $\rho^{[k]}$}\label{S:4}
From now on in the present paper, our goal is to prove Theorem \ref{T:specval} which is a special but more explicit case of Theorem \ref{T:1}. To prove our theorem, we first determine a certain subdomain of $\mathbb{C}_{\infty}^{rk+1}$ that lies in the domain of convergence of the logarithm function $\Log_{\rho^{[k]}}$ of $\rho^{[k]}$.

Let $\phi$ be the Drinfeld module of rank $r$ defined over $\mathbb{C}_{\infty}$ given by 
\begin{equation}\label{E:Drinfeldlogsec}
\phi_{\theta}=\theta+k_1\tau+\cdots+k_r\tau^r,\ \ \ k_r\not = 0
\end{equation}
and $\rho^{[k]}$ be the abelian Anderson module introduced in \S\ref{S:drinfeldtens}.  Let $q^{R_{\phi}}$ be the radius of convergence of $\log_{\phi}$. We also define $\mathfrak{r}_r:=R_{\phi}$ and for $1\leq m \leq r-1$, we set
 \begin{multline*}
 \mathfrak{r}_{r-m}:=\min\{ q^{m-1}(v_{\infty}(k_{r-m+1})+R_{\phi}), q^{m-2}(v_{\infty}(k_{r-m+2})+R_{\phi})\\
 ,\dots,q(v_{\infty}(k_{r-1})+R_{\phi}),v_{\infty}(k_r)+R_{\phi}\}.
 \end{multline*}

The following proposition, generalizing \cite[Prop. 2.4.1]{AT90} to Drinfeld modules of arbitrary rank tensored with the Carlitz tensor powers, will be proved in this section.
\begin{proposition}[{cf. \cite[Prop. 2.4.1]{AT90}}] \label{P:roclog} Let $\omega=(\omega_1,\dots,\omega_{rk+1})^{\tr}\in \mathbb{C}_{\infty}^{rk+1}$ be such that
\[
\log_q|\omega_1|<R_{\phi}+\frac{k}{q-1} 
\]
and
\[
\log_q|\omega_{ru+j+1}|<q^{j-r}\mathfrak{r}_{j}+ u+1+\frac{k}{q-1}, \ \ 0\leq u \leq k-1, \ \ 1\leq j\leq r.
\]
Then $\Log_{\rho^{[k]}}(\omega)$ converges on $\mathbb{C}_{\infty}^{rk+1}$. 
\end{proposition}

We start with an analysis done in \cite[\S6.1]{GN24}. Let $H_k:=\Mat_{1\times (rk+1)}(\mathbb{C}_{\infty}[\sigma])$ be the $t$-comotive of $\rho^{[k]}$ defined as in \S\ref{S:dualmot}. For $1\leq i \leq rk+1$, let $s_{i}\in \Mat_{1\times (rk+1)}(\mathbb{F}_q)$ be the $i$-th unit vector. Then note that the set $\{s_{1},\dots,s_{rk+1}\}$ forms a $\mathbb{C}_{\infty}[\sigma]$-basis for $H_k$. We consider the matrix $B\in \GL_r(\mathbb{C}_{\infty})$ whose inverse is given by
\begin{equation}\label{D:B def}
B^{-1} = \begin{pmatrix}
0 & k_2^{(-1)}&k_3^{(-2)} &\dots & k_{r-1}^{(-r+2)} & k_r^{(-r+1)}\\
0 & k_3^{(-1)}&k_4^{(-2)} & \dots& k_r^{(-r+2)}&\\
\vdots&\vdots&&\reflectbox{$\ddots$}&&\\
\vdots&\vdots&\reflectbox{$\ddots$}&&&\\
0 & k_r^{(-1)} &&&\\
1 &  & &&&
\end{pmatrix}.
\end{equation}
We further define 
\[
\mathfrak{I}:=(B^{-1})^{(-1)}(t-\theta)^k\Phi B. 
\]
For each $1\leq j \leq r$, set $\mathfrak{h}_j:=s_{r(k-1)+j+1}$ and define the $\mathbb{C}_{\infty}[t]$-basis $\mathfrak{h}$ given by $ \mathfrak{h}:=[s_{r(k-1)+2}, \dots, s_{rk+1}]^{\tr}$. Then we have 
\[
\sigma \mathfrak{h}=\mathfrak{I}\mathfrak{h}.
\]
We consider the map
\[
\iota: \Mat_{1\times r}(\CC_{\infty}[t]) \rightarrow  \Mat_{1\times (rk+1)}(\CC_{\infty}[\sigma])
\]
by setting for $[\mathfrak{g}_1, \dots, \mathfrak{g}_r] \in \Mat_{1\times r}(\mathbb{C}_{\infty}[t])$,
\[
\iota([\mathfrak{g}_1, \dots, \mathfrak{g}_r]) := \mathfrak{g}_1 \cdot \mathfrak{h}_1+ \dots+\mathfrak{g}_r \cdot \mathfrak{h}_r.
\]
We also consider the $\CC_{\infty}$-linear map 
\[\delta_0:\Mat_{1\times (rk+1)}(\CC_{\infty}[\sigma])\to \Mat_{(rk+1)\times 1}(\CC_{\infty})
\]
given by 
\[
\delta_0\left(\sum_{i\geq 0}B_i\sigma^i\right) := B_0^{\tr},\ \ B_i\in \Mat_{1\times rk+1}(\CC_{\infty}), \ \ B_i=0 \text{ for } i\gg0 .
\]
Using the exact same idea in the proof of \cite[Prop.~5.2]{Gez21}, we have
\begin{equation}\label{E:deltaiota}
\delta_0\circ\iota([\mathfrak{g}_1, \dots, \mathfrak{g}_r]) = [\partial_t^k(\mathfrak{g}_r),\partial_t^{k-1}(\mathfrak{g}_1), \dots, \partial_t^{k-1}(\mathfrak{g}_r),\dots, \dots, \partial_t^1(\mathfrak{g}_1), \dots,  \partial_t^1(\mathfrak{g}_r), \mathfrak{g}_1, \dots, \mathfrak{g}_r ]^{\tr} |_{t=\theta}.
\end{equation}
Moreover, there exists a unique extension of $\delta_0\circ\iota$ given by \[
\widehat{\delta_0\circ\iota}: \Mat_{1\times r}(\TT_{\theta}) \rightarrow \Mat_{(rk+1)\times 1}(\CC_{\infty}),
\]
where $\TT_{\theta}$ is the ring of rigid analytic functions on $\mathbb{C}_{\infty}$ convergent on the disk of radius $|\theta|=q$ centered at the origin.

\textit{A partition of a set $S$} is a disjoint union of sets whose union is equal to $S$. For any set $S\subset \mathbb{Z}$ and an integer $j$, we let $S+j:=\{s+j\ \ |s\in S\}$. For any $r\in \mathbb{Z}_{\geq 1}$, set $P_r(0):=\{(\emptyset,\dots,\emptyset)\}$ and for $n\in \mathbb{Z}_{\geq 1}$, we let $P_r(n)$ be the collection of sets $(S_1,\dots,S_r)$ so that each $S_i$ is a subset of $\{0,1,\dots,n-1\}$ and the tuple $\{S_i+j\ \ | 1\leq i \leq r, \ \ 0\leq j \leq i-1\}$ forms a partition for $\{0,1,\dots,n-1\}$. Following the notation in \cite[(6.4)]{EGP14}, we then let
\[
\beta_n(t) := \sum_{(S_1,\dots,S_r)\in P_r(n)} \prod_{i=1}^r \prod_{j\in S_i} \frac{k_i^{q^j}}{t-\theta^{q^{i+j}}}\in \mathbb{C}_{\infty}(t).
\]
For $n\geq 1$, let 
\[
R_n:=(\alpha^{[n]}_{\ell m}):=(\Phi^{-1})^{(n)}\cdots (\Phi^{-1})^{(1)}\in \GL_r(\mathbb{C}_{\infty}(t)).
\]
For later use, we state the following lemma which slightly generalizes and corrects couple of typos in \cite[Lem. 4.1]{GG25}. For completeness, we provide the details of its proof. 
\begin{lemma}For each $1\leq \ell,m \leq r$ and $n\geq r-1$, we have
  \begin{equation}\label{E:deseq}
	\alpha^{[n]}_{\ell m}=
\sum_{v=0}^{r-m}\frac{\beta^{(v+1)}_{n-\ell-v}(t)}{t-\theta^{q^{v+1}}}k_{m+v}^{(1-m)}.
	\end{equation} 
\end{lemma}
\begin{proof} For each $1\leq \ell \leq r$ and $n\geq r-1$, our first goal is to show the following: We have\begin{equation}\label{E:log3}
\alpha^{[n]}_{\ell1}=\beta_{n-(\ell-1)}(t),\ \ \ \ 
\alpha^{[n]}_{\ell r}=\frac{k_r^{(-r+1)}\beta_{n-\ell}^{(1)}(t)}{t-\theta^q}
\end{equation}
and for $1\leq m \leq r-2$ and $n\geq m+1$, 
\begin{equation}\label{E:log4}
	\alpha^{[n]}_{\ell (r-m)}=
\sum_{v=0}^{m}\frac{\beta^{(v+1)}_{n-\ell-v}(t)}{t-\theta^{q^{v+1}}}k_{r-m+v}^{(1-(r-m))}.
	\end{equation}
 The first equality in \eqref{E:log3} follows from \cite[Lem. 3.1.4]{Che24}. For the second equality,  when $n\geq 2$, we have 
	\begin{multline}\label{E:prod}
	R_{n}=R_{n-1}^{(1)}(\Phi^{-1})^{(1)}\\
    =\begin{pmatrix}
	\beta_{n-1}(t)&\alpha_{12}^{[n-1]}&\cdots& \cdots &\alpha_{1r}^{[n-1]}\\
	\beta_{n-2}(t)&\alpha_{22}^{[n-1]} & & &\alpha_{2r}^{[n-1]}\\
	\vdots&\vdots & & &\vdots \\
	\vdots&\vdots &  & &\vdots \\
	\beta_{n-1-(r-1)}(t)&\alpha_{r2}^{[n-1]}&\cdots &\cdots & \alpha_{rr}^{[n-1]}
	\end{pmatrix}^{(1)}\begin{pmatrix}
	\frac{k_1}{t-\theta^q}&\frac{k_2^{(-1)}}{t-\theta^q}&\cdots& \frac{k_{r-1}^{(-(r-2))}}{t-\theta^q}&\frac{k_r^{(-(r-1))}}{t-\theta^q}\\
	1& & & & \\
	 &\ddots & & & \\
	 & & \ddots & & \\
	 & & & 1& 0
	\end{pmatrix}.
	\end{multline}
	Now, the second equality in \eqref{E:log3} follows from the first equality and \eqref{E:prod}. We show the equality in \eqref{E:log4}. By the first equality in \eqref{E:log3} and \eqref{E:prod}, we have
    \[
    \alpha_{\ell(r-m)}^{[n]}=\frac{\mathcal{B}^{(1)}_{n-1-(\ell-1)}(t)}{t-\theta^q}k_{r-m}^{(1-(r-m))}+(\alpha_{\ell(r-(m-1))}^{[n-1]})^{(1)}.
    \]
    For $m=1$ and $n\geq 2$, we have 
    \begin{multline*}
\alpha_{\ell(r-1)}^{[n]}=\frac{\mathcal{B}^{(1)}_{n-1-(\ell-1)}(t)}{t-\theta^q}k_{r-1}^{(1-(r-1))}+(\alpha_{\ell r}^{[n-1]})^{(1)} \\
=\frac{\mathcal{B}^{(1)}_{n-\ell}(t)}{t-\theta^q}k_{r-1}^{(1-(r-1))}+\left( \frac{k_r^{(-(r-1))}\beta_{n-1-\ell}^{(1)}(t)}{t-\theta^q} \right)^{(1)}=\frac{\mathcal{B}^{(1)}_{n-\ell}(t)}{t-\theta^q}k_{r-1}^{(1-(r-1))}+\frac{\beta_{n-1-\ell}^{(2)}(t)}{t-\theta^{q^2}}k_r^{(1-(r-1))}\\
=\sum_{v=0}^{1}\frac{\beta^{(v+1)}_{n-\ell-v}(t)}{t-\theta^{q^{v+1}}}k_{r-1+v}^{(1-(r-1))}.
\end{multline*}
Now, for any $2\leq m \leq r-2$ and $n\geq m+1$, applying the above idea recursively, we obtain
\begin{multline*}
\alpha_{\ell(r-m)}^{[n]}=\frac{\mathcal{B}^{(1)}_{n-1-(\ell-1)}(t)}{t-\theta^q}k_{r-m}^{(1-(r-m))}+(\alpha_{\ell(r-(m-1))}^{[n-1]})^{(1)}\\
=\frac{\mathcal{B}^{(1)}_{n-\ell}(t)}{t-\theta^q}k_{r-m}^{(1-(r-m))}+\left(\sum_{v=0}^{m-1}\frac{\beta^{(v+1)}_{n-1-\ell-v}(t)}{t-\theta^{q^{v+1}}}k_{r-(m-1)+v}^{(1-(r-(m-1)))}\right)^{(1)}\\
=\frac{\mathcal{B}^{(1)}_{n-\ell}(t)}{t-\theta^q}k_{r-m}^{(1-(r-m))}+\sum_{v=0}^{m-1}\frac{\beta^{(v+2)}_{n-1-\ell-v}(t)}{t-\theta^{q^{v+2}}}k_{r-(m-1)+v}^{(1-(r-m))}=\sum_{v=0}^{m}\frac{\beta^{(v+1)}_{n-\ell-v}(t)}{t-\theta^{q^{v+1}}}k_{r-m+v}^{(1-(r-m))}\end{multline*}
as desired. Now, after changing the index of the summation in \eqref{E:log4}, when $2\leq m \leq r-1$, \eqref{E:deseq}  follows from \eqref{E:log4} and when $m=r$, \eqref{E:deseq} follows from the second equality in \eqref{E:log3}. Finally, when $m=1$, \eqref{E:deseq} follows from the first equality in \eqref{E:log3} and \cite[Lem. 6.12(a)]{EGP14}.
\end{proof}

Using \eqref{E:log3} and \eqref{E:log4}, observe that 
\begin{multline}\label{E:matrixmult}
R_n=\begin{pmatrix} \beta_n(t)&\frac{\beta_{n-(r-1)}^{(r-1)}(t)}{t-\theta^{q^{r-1}}}&\cdots & \frac{\beta_{n-1}^{(1)}(t)}{t-\theta^{q}}\\
\beta_{n-1}(t)&\frac{\beta_{n-r}^{(r-1)}(t)}{t-\theta^{q^{r-1}}}&\cdots & \frac{\beta_{n-2}^{(1)}(t)}{t-\theta^{q}}\\
\vdots&\vdots & & \vdots\\
\beta_{n-(r-1)}(t)&\frac{\beta_{n-2(r-1)}^{(r-1)}(t)}{t-\theta^{q^{r-1}}}&\cdots & \frac{\beta_{n-r}^{(1)}(t)}{t-\theta^{q}}
\end{pmatrix}\begin{pmatrix}
    1& & & & \\
     &k_r^{(-1)} & & & \\
     &k_{r-1}^{(-1)} & \ddots & & \\
 & \vdots& \ddots &\ddots & \\
 & k_2^{(-1)} &\cdots & k_{r-1}^{(-r+2)}&k_r^{(-r+1)}
\end{pmatrix}\\
=\begin{pmatrix} \beta_n(t)&\beta_{n-(r-1)}^{(r-1)}(t)&\cdots & \beta_{n-1}^{(1)}(t)\\
\beta_{n-1}(t)&\beta_{n-r}^{(r-1)}(t)&\cdots & \beta_{n-2}^{(1)}(t)\\
\vdots&\vdots & & \vdots\\
\beta_{n-(r-1)}(t)&\beta_{n-2(r-1)}^{(r-1)}(t)&\cdots & \beta_{n-r}^{(1)}(t)
\end{pmatrix}\begin{pmatrix}
    1& & & & \\
     & \frac{1}{t-\theta^{q^{r-1}}} & & & \\
     & & \ddots & & \\
     & & & & \frac{1}{t-\theta^q}
\end{pmatrix}  \begin{pmatrix}
    & & & 1\\
    & & \iddots & \\
    & \iddots & & \\
    1 & & & 
\end{pmatrix}B^{-1}\\
=\begin{pmatrix} \beta_n(t)&\beta_{n-(r-1)}^{(r-1)}(t)&\cdots & \beta_{n-1}^{(1)}(t)\\
\beta_{n-1}(t)&\beta_{n-r}^{(r-1)}(t)&\cdots & \beta_{n-2}^{(1)}(t)\\
\vdots&\vdots & & \vdots\\
\beta_{n-(r-1)}(t)&\beta_{n-2(r-1)}^{(r-1)}(t)&\cdots & \beta_{n-r}^{(1)}(t)
\end{pmatrix}\begin{pmatrix}
    & & & 1\\
    & & \frac{1}{t-\theta^{q^{r-1}}} & \\
    & \iddots & & \\
    \frac{1}{t-\theta^q} & & & 
\end{pmatrix}B^{-1},
\end{multline}
where we remind the reader of the definition of $B^{-1}$ from \eqref{D:B def}.

For $1\leq j \leq r$, recall that $\mathfrak{e}_j\in \Mat_{r\times 1}(\mathbb{F}_q)$ is the $j$-th unit vector. For any $n\geq 1$, let $\mathfrak{I}^{[n]}:=(\mathfrak{I}^{-1})^{(n)}\cdots (\mathfrak{I}^{-1})^{(1)}\in \GL_r(\mathbb{C}_{\infty}(t))$. Using \cite[(43)]{GN24},  \cite[Cor. 4.5]{Gre24} and letting $L_{\ell,n}\in \Mat_{(rk+1)\times 1}(\mathbb{C}_{\infty})$ be the $\ell$-th column of the $n$-th coefficient $L_n$ of $\Log_{\rho^{[k]}}$, we  have 
\begin{equation}\label{E:log1}
L_{1,n}=\widehat{\delta_0\circ\iota}(\mathfrak{e}_{r}(t-\theta^{q^n})^{k}\mathfrak{I}^{[n]})
\end{equation}
and 
\begin{equation}\label{E:log2}
L_{ru+j+1,n}=\widehat{\delta_0\circ\iota}(\mathfrak{e}_{j}(t-\theta^{q^n})^{k-(u+1)}\mathfrak{I}^{[n]}).
\end{equation}

To prove Proposition \ref{P:roclog}, we finally need our next lemma on the size of the hyperderivatives of $\beta_n(t)$ evaluated at $t=\theta$.
\begin{lemma}\label{L:hypebeta}  For any $m\geq 0$, we have 
\[
\log_q|\partial_t^m(\beta_n(t))|_{t=\theta}| \leq -q^{n}R_{\phi}+\mathfrak{R},
\]
where $\mathfrak{R}$ is a real number independent of $n$.
\end{lemma}
\begin{proof} First, we note that
\begin{equation}\label{E:thypder}
\partial_t^m\left(\frac{1}{t-\theta^{q^{\ell}}}\right)=\frac{(-1)^m}{(t-\theta^{q^{\ell}})^{m+1}}, \ \ m\geq 0, \ \ \ell\geq 0.
\end{equation}
For each $1\leq i \leq r$ and $(S_1,\dots,S_r)\in P_r(n)$, let $S_i=\{j_{i1},\dots,j_{i\mathfrak{s}_i}\}$. When $m=0$, the lemma follows from \cite[Lem. 6.7(b)]{EGP14} and \cite[Cor. 4.5]{KP23}. Let $m\geq 1$. By the Leibniz rule \eqref{E:Leibniz} and \eqref{E:thypder}, we have 
\begin{multline*}
\partial_t^m(\beta_n(t))=\partial_t^m\left(\sum_{(S_1,\dots,S_r)\in P_r(n)} \prod_{i=1}^r \prod_{j\in S_i} \frac{k_i^{q^j}}{t-\theta^{q^{i+j}}}\right)=\sum_{(S_1,\dots,S_r)\in P_r(n)}\partial_t^m\left( \prod_{i=1}^r \prod_{j\in S_i} \frac{k_i^{q^j}}{t-\theta^{q^{i+j}}}\right)\\
=\sum_{(S_1,\dots,S_r)\in P_r(n)} \sum_{\substack{\ell_{11},\dots,\ell_{1\mathfrak{s}_1},\dots, \ell_{r1},\dots,\ell_{r\mathfrak{s}_r} \geq 0\\ \ell_{11}+\cdots+\ell_{1\mathfrak{s}_1}+\cdots+\ell_{r1}+\cdots+\ell_{r\mathfrak{s}_r}=m}}
\partial_t^{\ell_{11}}\left(\frac{k_1^{q^{j_{11}}}}{t-\theta^{q^{1+j_{11}}}}\right)\
\cdots \partial_t^{\ell_{1\mathfrak{s}_1}}\left(\frac{k_1^{q^{j_{1\mathfrak{s}_1}}}}{t-\theta^{q^{1+j_{1\mathfrak{s}_1}}}}\right) \times \\ \cdots \partial_t^{\ell_{r1}}\left(\frac{k_r^{q^{j_{r1}}}}{t-\theta^{q^{r+j_{r1}}}}\right)\cdots \partial_t^{\ell_{rr}}\left(\frac{k_r^{q^{j_{r\mathfrak{s}_r}}}}{t-\theta^{q^{r+j_{r\mathfrak{s}_r}}}}\right)\\
=\sum_{(S_1,\dots,S_r)\in P_r(n)}\prod_{i=1}^r \prod_{j\in S_i}\frac{k_i^{q^j}}{t-\theta^{q^{i+j}}} \sum_{\substack{\ell_{11},\dots,\ell_{1\mathfrak{s}_1},\dots, \ell_{r1},\dots,\ell_{r\mathfrak{s}_r} \geq 1\\ \ell_{11}+\cdots+\ell_{1\mathfrak{s}_1}+\cdots+\ell_{r1}+\cdots+\ell_{r\mathfrak{s}_r}=m}} \prod_{\substack{1\leq i \leq r\\ 1\leq \nu \leq \mathfrak{s}_i}}   \frac{(-1)^{\ell_{i\nu}}}{(t-\theta^{q^{i+j_{i\nu}}})^{\ell_{i\nu}}}.\end{multline*}As $|\theta-\theta^{i+j_{i\nu}}|^{-1}\leq 1$ for $1\leq i \leq r$ and $1\leq \nu\leq \max\{\mathfrak{s}_1,\dots,\mathfrak{s}_r\}$, above calculation, together with \cite[Lem. 6.7(b)]{EGP14} and \cite[Cor. 4.5]{KP23}, immediately implies the lemma.
\end{proof}

\begin{proof}[{Proof of Proposition \ref{P:roclog}}] Our goal, by using \eqref{E:log1} and \eqref{E:log2}, is to estimate the norm of each entry of $L_n$. For each $n\geq 1$, recall the polynomials $\mathbb{S}_n\in \mathbb{C}_{\infty}[t]$ defined in \eqref{E:sn}. Observe, by using \eqref{E:matrixmult}, that 
\begin{multline}\label{E:invmat}
\mathfrak{I}^{[n]}=(\mathfrak{I}^{-1})^{(n)}\cdots (\mathfrak{I}^{-1})^{(1)}=\mathbb{S}_n^{-k}(B^{-1})^{(n)}(\Phi^{-1})^{(n)}\cdots (\Phi^{-1})^{(1)}B=\mathbb{S}_n^{-k}(B^{-1})^{(n)}R_nB\\
=\mathbb{S}_n^{-k}(B^{-1})^{(n)}\begin{pmatrix} \beta_n(t)&\beta_{n-(r-1)}^{(r-1)}(t)&\cdots & \beta_{n-1}^{(1)}(t)\\
\beta_{n-1}(t)&\beta_{n-r}^{(r-1)}(t)&\cdots & \beta_{n-2}^{(1)}(t)\\
\vdots&\vdots & & \vdots\\
\beta_{n-(r-1)}(t)&\beta_{n-2(r-1)}^{(r-1)}(t)&\cdots & \beta_{n-r}^{(1)}(t)
\end{pmatrix}\begin{pmatrix}
    & & & 1\\
    & & \frac{1}{t-\theta^{q^{r-1}}} & \\
    & \iddots & & \\
    \frac{1}{t-\theta^q} & & & 
\end{pmatrix}\\
=\mathbb{S}_n^{-k}\begin{pmatrix}
0 & k_2^{(n-1)}&k_3^{(n-2)} &\dots & k_{r-1}^{(n-r+2)} & k_r^{(n-r+1)}\\
0 & k_3^{(n-1)}&k_4^{(n-2)} & \dots& k_r^{(n-r+2)}&\\
\vdots&\vdots&&\reflectbox{$\ddots$}&&\\
\vdots&\vdots&\reflectbox{$\ddots$}&&&\\
0 & k_r^{(n-1)} &&&\\
1 &  & &&&
\end{pmatrix}\begin{pmatrix} \frac{\beta_{n-1}^{(1)}(t)}{t-\theta^{q}}&\cdots &\frac{\beta_{n-(r-1)}^{(r-1)}(t)}{t-\theta^{q^{r-1}}} & \beta_n(t)\\
\frac{\beta_{n-2}^{(1)}(t)}{t-\theta^{q}}&\cdots&\frac{\beta_{n-r}^{(r-1)}(t)}{t-\theta^{q^{r-1}}} &\beta_{n-1}(t)\\
\vdots& & \vdots& \vdots\\
\frac{\beta_{n-r}^{(1)}(t)}{t-\theta^{q}}&\cdots & \frac{\beta_{n-2(r-1)}^{(r-1)}(t)}{t-\theta^{q^{r-1}}}&\beta_{n-(r-1)}(t)
\end{pmatrix}\\
=\mathbb{S}_n^{-k}
\begin{pmatrix}
\frac{\sum_{i=2}^{r}k_i^{(n-i+1)}\beta_{n-i}^{(1)}(t)}{t-\theta^q}& & \cdots & \frac{\sum_{i=2}^{r}k_i^{(n-i+1)}\beta_{n-i-r+2}^{(r-1)}(t)}{t-\theta^{q^{r-1}}}&\sum_{i=2}^{r}k_i^{(n-i+1)}\beta_{n-i+1}(t) 
\\
\frac{\sum_{i=3}^{r}k_i^{(n-i+2)}\beta_{n-i+1}^{(1)}(t)}{t-\theta^q}& & \cdots & \frac{\sum_{i=3}^{r}k_i^{(n-i+2)}\beta_{n-i-r+3}^{(r-1)}(t)}{t-\theta^{q^{r-1}}}&\sum_{i=3}^{r}k_i^{(n-i+2)}\beta_{n-i+2}(t) \\
\vdots & &  & \vdots & \vdots\\
\frac{k_r^{(n-1)}\beta_{n-2}^{(1)}(t)}{t-\theta^q} &  &\cdots  &  \frac{k_r^{(n-1)}\beta_{n-r}^{(r-1)}(t)}{t-\theta^{q^{r-1}}} & k_r^{(n-1)}\beta_{n-1}(t)\\
\frac{\beta_{n-1}^{(1)}(t)}{t-\theta^q}& &\cdots & \frac{\beta_{n-(r-1)}^{(r-1)}(t)}{t-\theta^{q^{r-1}}}& \beta_n(t)
\end{pmatrix}.
\end{multline}
Moreover, letting $|\mathfrak{I}^{[n]}_{j}|$ be the maximum norm among the elements in the $j$-th row of $\mathfrak{I}^{[n]}$ as well as using Lemma \ref{L:hypebeta} and  \eqref{E:invmat}, we further obtain 
\begin{equation}\label{E:log7}
\log_q|\mathfrak{I}^{[n]}_{j}|\leq -q^{n-(r-j)}\mathfrak{r}_j-\frac{kq^{n+1}}{q-1}+\mathfrak{O}_j
\end{equation}
where $\mathfrak{O}_j$ is a real number independent of $n$. Thus, 
the desired result follows from \eqref{E:log1}, \eqref{E:log2} and \eqref{E:log7}.
\end{proof}


\subsection{Drinfeld modules with $\mathfrak{m}_{0}=0$}\label{S:4.1}In what follows, we  assume that $\phi$ is a  Drinfeld module over $\mathbb{C}_{\infty}$ given as in \eqref{E:Drinfeldlogsec} so that  $\mathfrak{m}_{0}=0$. In other words, letting $\lambda_1,\dots,\lambda_r$ be generators of the period lattice of $\phi$ as an $\bA$-module, we assume that $\lambda_{i}/\theta$ lies in the domain of convergence of $\log_{\phi}$ for all $1\leq i\leq r$.

Our goal in the following subsections is to derive bounds for the periods of $\rho^{[k]}$ described in \eqref{E:periodsofext} where $\phi$ is a particular Drinfeld module with $\mathfrak{m}_{0}=0$.

We start with making some analysis on the Anderson-Thakur function. Note, by \cite[Prop. 3.2]{EGP14}, that for $\nu \geq 1$,  one can express $\omega_C^{\nu}$ as an infinite sum given by 
\begin{equation}\label{E:omegaexp}
\omega_C^{\nu}=\sum_{n=-\nu}^{\infty}c_{n,\nu}(t-\theta)^n, \ \ c_{n,\nu}\in \mathbb{C}_{\infty}.
\end{equation}
In our next lemma, we  estimate the norm of $c_{n,\nu}$ when $n\leq 0$.
\begin{lemma}\label{L:estimates} For $\nu\geq 1$, let
\[
\omega_C^{\nu}=\frac{c_{-\nu,\nu}}{(t-\theta)^{\nu}}+\cdots +\frac{c_{-1,\nu}}{t-\theta}+R_{\nu}(t)
\]
where $R_{\nu}(t):=\sum_{n=0}^{\infty}c_{n,\nu}(t-\theta)^{n}$. Then for $0\leq \ell \leq \nu$, we have
\[
\log_q|c_{-\ell,\nu}|\leq \ell-\nu+\frac{\nu q}{q-1}.
\] 
\end{lemma}
\begin{proof} By Lemma \ref{L:AGF}(i), note that
\begin{equation}\label{E:ATF}
\omega_C=\sum_{n=0}^{\infty}\Exp_{C}\left(\frac{\tilde{\pi}}{\theta^{n+1}}\right)t^{n}=-\frac{\tilde{\pi}}{t-\theta}+\sum_{i=1}^{\infty}\frac{\alpha_i\tilde{\pi}^{q^i}}{\theta^{q^i}-t}
=\sum_{n=-1}^{\infty}c_{n,1}(t-\theta)^n\end{equation}
so that $c_{-1,1}=-\tilde{\pi}$ and, by \cite[Prop. 3.2.3]{CGM19}, for $n\geq 0$, we have
\[
c_{n,1}=\partial^{n+1}_t((t-\theta)\omega_C)|_{t=\theta}=\partial_t^{n}(\omega_C)|_{t=\theta}=\left(\sum_{i=n}^{\infty}\Exp_{C}\left(\frac{\tilde{\pi}}{\theta^{i+1}}\right)\binom{i}{n}t^{i-n}\right)|_{t=\theta}.
\]
Since $|\tilde{\pi}|=q^{q/(q-1)}$, we have 
\begin{equation}\label{E:ineq0}
\log_q|c_{-1,1}|=\frac{q}{q-1}.
\end{equation}
On the other hand, by \cite[Prop./Def. 2.43]{AT90}, we see that $\tilde{\pi}/\theta^n$ for $n\geq 1$ lies in the domain of convergence of $\Log_{C}$. Since, on the domain of convergence of $\Log_C$, $\Exp_{C}$ is an isometric map (see for example \cite[Prop. 4.14.12]{Gos96}), we have $|\Exp_C(\tilde{\pi}/\theta^n)|=|\tilde{\pi}/\theta^n|$ for all $n\geq 1$ and hence
\begin{equation}\label{E:ineq1}
\log_q|c_{n,1}|\leq \max\left\{\left|\frac{\tilde{\pi}}{\theta^{i+1}}\right||\theta|^{i-n} \ \ | \ \ i\geq n\right\}= \left|\frac{\tilde{\pi}}{\theta^{n+1}}\right|\leq \frac{q}{q-1}-(n+1).
\end{equation}
Using \eqref{E:ineq0} and \eqref{E:ineq1},  we now observe that 
\begin{multline*}
|c_{-\ell,\nu}|\\
\leq \max\{|c_{-1,1}^{m_0}c_{\ell_1,1}^{m_1}\cdots c_{\ell_s,1}^{m_s}| \ \ | \ \ -m_0+\sum_{\mu=1}^s\ell_{\mu}m_{\mu}=-\ell, \ \ \sum_{\mu=0}^sm_{\mu}=\nu, \ \ m_0,\dots,m_s\geq 0,\ \ s\geq 0 \}\\
\leq \frac{q}{q-1}(m_0+\cdots+m_s)-((-1+1)m_0+(\ell_1+1)m_1+\cdots+(\ell_s+1)m_s)\\
\leq \frac{\nu q}{q-1} +m_0-(\ell_1m_1+\cdots+\ell_sm_s)-(m_0+m_1+\cdots+m_s)=\ell-\nu+\frac{\nu q}{q-1}
\end{multline*}
as desired.
\end{proof}

For each $1\leq j \leq r$, $n\geq 0$ and $0\leq \ell \leq r-1$, we define $\mathfrak{C}^{[\ell]}_{n,j}$ so that $\mathfrak{C}_{-1,j}^{[0]}:=-\lambda_j$, 
\[
f_{\lambda_j}(t)=-\frac{\lambda_j}{t-\theta}+\sum_{n=0}^{\infty}\mathfrak{C}^{[0]}_{n,j}(t-\theta)^{n}=\sum_{n=-1}^{\infty}\mathfrak{C}^{[0]}_{n,j}(t-\theta)^{n}\]
and 
\[
f^{(\ell)}_{\lambda_j}(t)=\sum_{n=0}^{\infty}\mathfrak{C}^{[\ell]}_{n,j}(t-\theta)^{n}.
\]
Note, by \cite[Prop. 3.2.3]{CGM19} and the Leibniz rule, that $
\mathfrak{C}_{n,j}^{[0]}=\partial_t^{n+1}((t-\theta)f_{\lambda_j}(t))|_{t=\theta}=\partial_t^{n}(f_{\lambda_j}(t))|_{t=\theta}$
and $
\mathfrak{C}_{n,j}^{[\ell]}=\partial_t^n(f_{\lambda_j}^{(\ell)}(t))|_{t=\theta}.
$
\begin{lemma}\label{L:lemmalog1}For $1\leq j \leq r$, we have $||f_{\lambda_j}(t)||=\left|\frac{\lambda_j}{\theta}\right| $.
Moreover, for $1\leq \ell \leq r-1$ and $n\geq 0$, we have 
\[
|\mathfrak{C}_{n,j}^{[\ell]}|\leq \left|\frac{\lambda_{j}}{\theta^{n+1}}\right|^{q^{\ell}}=\left|\frac{\lambda_{j}}{\theta}\right|^{q^{\ell}}\frac{1}{|\theta|^{nq^{\ell}}}.
\]
\end{lemma}
\begin{proof} 
    Note that, since $\mathfrak{m}_{0}=0$ and $\Exp_{\phi}$ is an isometric map on the domain of convergence of $\Log_{\phi}$, for $i\geq 0$, we have
    \begin{equation}\label{E:normofexp}
   \left| \Exp_{\phi}\left(\frac{\lambda_j}{\theta^{i+1}}\right)\right|=\left|\frac{\lambda_j}{\theta^{i+1}}\right|.
    \end{equation}
Hence, we obtain $||f_{\lambda_j}(t)||= \left|\frac{\lambda_j}{\theta}\right| $. We now prove the second assertion. 
First, note that 
    \[
    \mathfrak{C}_{n,j}^{[0]}=\partial_t^{n}(f_{\lambda_j}(t))|_{t=\theta}=\left(\sum_{i=n}^{\infty}\Exp_{\phi}\left(\frac{\lambda_j}{\theta^{i+1}}\right)\binom{i}{n}t^{i-n}\right)|_{t=\theta}.
    \]
    Thus, by \eqref{E:normofexp}, we obtain
    \[
|\mathfrak{C}_{n,j}^{[0]}|\leq \max\left\{\left|\frac{\lambda_j}{\theta^{i+1}}\right||\theta|^{i-n} \ \ | \ \ i\geq n\right\}= \left|\frac{\lambda_{j}}{\theta^{n+1}}\right|.
\]
Similarly, we have 
\[
\mathfrak{C}_{n,j}^{[\ell]}=\partial_t^{n}(f_{\lambda_j}^{(\ell)}(t))|_{t=\theta}=\left(\sum_{i=n}^{\infty}\Exp_{\phi}\left(\frac{\lambda_j}{\theta^{i+1}}\right)^{q^{\ell}}\binom{i}{n}t^{i-n}\right)|_{t=\theta}.
\]
Using \eqref{E:normofexp}, it implies 
\[
|\mathfrak{C}_{n,j}^{[\ell]}|\leq \max\left\{\left|\frac{\lambda_j}{\theta^{i+1}}\right|^{q^{\ell}}|\theta|^{i-n} \ \ | \ \ i\geq n\right\}= \left|\frac{\lambda_{j}}{\theta^{n+1}}\right|^{q^{\ell}}
\]
as desired.
\end{proof}

We  state estimates on the size of the generators of the period lattice   $\Lambda_{\rho^{[k]}}$. For each $1\leq j \leq r$, recall the period $\varpi_{k,j}$ of $\Lambda_{\rho^{[k]}}$ given in \eqref{E:periodsofext}.

\begin{lemma}\label{L:ourperiod} Let $\varpi_{k,j}=(\mathfrak{z}_1,\dots,\mathfrak{z}_{rk+1})^{\tr}\in \mathbb{C}_{\infty}^{rk+1}$ and set $\mathfrak{z}_{n}:=0$ if $n>rk+1$. Then 
\[
\log_q|\mathfrak{z}_{\mu+ru}|\leq q^{\mu-1}(\log_q|\lambda_j|-1)+u+1+\frac{k}{q-1},\ \ 0\leq u \leq k, \ \ 1\leq \mu \leq r.
\]
\end{lemma}
\begin{proof}
By \eqref{E:periodsofext}, Lemma \ref{L:estimates} and Lemma  \ref{L:lemmalog1}, we obtain
    \begin{multline*}
    \log_q|\mathfrak{z}_{1+ru}|=\log_q|\Res_{t=\theta}((t-\theta)^{u}\omega_C^kf_j(t)dt)| \\\leq \max\left\{ |c_{-m,k}\mathfrak{C}_{n,j}^{[0]}| \ \ | \ \ u\leq m 
    \leq k, \ \ n-m+u=-1 \right\}\\
    \leq 
    m-k+\frac{kq}{q-1}+(\log_q|\lambda_j|-1)-n=(R_{\phi}-1)+u+1+\frac{k}{q-1}.
    \end{multline*}
    Furthermore, for $1\leq u \leq k$,  we have
    \begin{multline*}
    \log_q|\mathfrak{z}_{\mu+ru}|=\log_q|\Res_{t=\theta}((t-\theta)^{u}\omega_C^kf_j^{(\mu-1)}(t)dt)| \\ \leq \max\left\{  |c_{-m,k}\mathfrak{C}^{[\mu-1]}_{n,j}| \ \ |  \ \  u\leq m
    \leq k, \ \ n-m+u=-1  \right\}\\
    \leq 
    m-k+\frac{kq}{q-1}+q^{\mu-1}(\log_q|\lambda_j|-1)-n=q^{\mu-1}(\log_q|\lambda_j|-1)+u+1+\frac{k}{q-1},
\end{multline*}
finishing the proof of the lemma.
    \end{proof}

\subsection{A certain class of Drinfeld modules with $\mathfrak{m}_{0}=0$}\label{S:4.2}
Let $\phi$ be a Drinfeld module of rank $r$  given by 
\begin{equation}\label{E:ourdm}
\phi_{t}=\theta+k_1\tau+\cdots+k_r\tau^r
\end{equation}
so that $|k_1|,\dots,|k_{r-1}|\leq 1$ and $k_r\in \mathbb{F}_q^{\times}$. 
Recall that $q^{R_{\phi}}$ is the radius of convergence of $\log_{\phi}$. In this case, we have
 \[
 R_{\phi}=\frac{q^r}{q^r-1}.
 \]
Recall the generators $\lambda_1,\dots,\lambda_r$ for the period lattice $\Lambda_{\phi}$ of $\phi$ over $A$. For any $1\leq j \leq r$, using the Newton polygon of $\phi_t$, we see that 
\begin{equation}\label{E:periodnorm}
\log_q|\lambda_j|=R_{\phi}.
\end{equation}
This implies that $\mathfrak{m}_{0}=0$. Furthermore, recalling the notation used in Proposition \ref{P:roclog},  we have 
\[
\mathfrak{r}_1=\cdots=\mathfrak{r}_r=R_{\phi}.
\]

For $k\geq 1$, let $\rho^{[k]}$ be the Anderson module defined in \S\ref{S:drinfeldtens} by using $\phi$ as in \eqref{E:ourdm}. Recall the set $\{\varpi_{k,1},\dots,\varpi_{k,rk+1}\}$ from Proposition \ref{P:periodbasis} which forms an $\bA$-module basis for $\Lambda_{\rho^{[k]}}$. 
\begin{corollary}\label{C:ourperiod} For each $n\geq 1$ and $1\leq j \leq r$, $\dmat_k[t]^{-n}\varpi_{k,j}$ lies in the domain of the convergence of $\Log_{\rho^{[k]}}$. 
\end{corollary}
\begin{proof} For each $n,\ell \geq 1$,let $\mathfrak{p}_{\ell,n}:=\partial_{t}^{\ell}(t^{-n})$ and observe that 
\[
\mathfrak{p}_{\ell,n}=\partial_{t}^{\ell}(t^{-n})=\sum_{m=1}^{\ell}\binom{\ell+1}{m+1}\frac{(-1)^m}{t^{n(m+1)}}\partial^{\ell}_{t}(t^{nm}).
\]

Moreover, note that if we set 
\[
\dmat_k[t]^{-n}\varpi_j= (
\mathfrak{u}_{1,n},\dots,
\mathfrak{u}_{rk+1,n})^{\tr},
\]
then 
\[
\mathfrak{u}_{1,n}=\sum_{\ell=0}^{k}\mathfrak{z}_{r\ell+1}\mathfrak{p}_{\ell,n}|_{t=\theta}
\]
and 
\[
\mathfrak{u}_{ru+j+1,n}=\sum_{\ell=u}^{k-1}\mathfrak{z}_{r\ell +j+1}\mathfrak{p}_{\ell-u,n}|_{t=\theta}, \ \ 1\leq j \leq r, \ \ 0\leq u \leq k-1.
\]
On the other hand, since $|\mathfrak{p}_{\ell,n}|\leq |\theta|^{-n-\ell}$ and $R_{\phi}=q^r/(q^r-1)$, by Lemma \ref{L:ourperiod}, we have
\begin{equation}\label{E:u1}
\log_q|\mathfrak{u}_{1,n}|\leq  R_{\phi}+\frac{k}{q-1}-n 
\end{equation}
and 
\begin{equation}\label{E:u2}
\log_q|\mathfrak{u}_{ru+j+1,n}|\leq q^{j-r}R_{\phi}+u+1+\frac{k}{q-1}-n.
\end{equation}
Thus, by Proposition \ref{P:roclog}, $\Log_{\rho^{[k]}}(\dmat_k[t]^{-n}\varpi_{k,j})$ converges.
\end{proof}

We finish this subsection by determining the integer $\mathfrak{m}_k$ corresponding to $\rho^{[k]}$. 
\begin{proposition}\label{P:norms} The following statements hold.
\begin{itemize}
    \item[(i)] Let $L_{\ell,m}$ be the $\ell$-th column of the $m$-th logarithm coefficient $L_{m}$ of $\rho^{[k]}$. Then 
    \[
    \log_q|L_{1,m}|\leq -q^m\left(R_{\phi}+\frac{k}{q-1}\right)+\mathfrak{P}_1
    \]
     and 
    \[
    \log_q|L_{ru+j+1,m}|\leq -q^m\left(q^{j-r}R_{\phi}+u+1+\frac{k}{q-1}\right)+\mathfrak{P}_{ru+j+1}, \ \ \  \ \  1\leq j \leq r, \ \ 0\leq u \leq k-1
    \]
    where $\mathfrak{P}_{1}$ and $\mathfrak{P}_{ru+j+1}$ are real numbers independent of $n$
    \item[(ii)] Let $\mathfrak{U}\subset \mathbb{C}_{\infty}^{rk+1}$ be the set of all $\mathfrak{z}=(z_1,\dots,z_{rk+1})\in \mathbb{C}_{\infty}^{rk+1}$ such that 
 \[
\log_q|z_1|\leq  R_{\phi}+\frac{k}{q-1}-1 
\]
and 
\[
\log_q|z_{ru+j+1}|\leq q^{j-r}R_{\phi}+u+\frac{k}{q-1} \ \ \  \ \  1\leq j \leq r, \ \ 0\leq u \leq k-1.
\]
Then $\Log_{\rho^{[k]}}(\mathfrak{U})\subseteq \mathfrak{U}$. 
Moreover, for any $1\leq j \leq r$ and $n\geq 1$,  we have 
    \[
    \left|\Log_{\rho^{[k]}}(d_k[t]^{-n}\varpi_{k,j})\right|=\left|d_k[t]^{-n}\varpi_{k,j}\right|.
    \]
    Furthermore,    $\mathfrak{m}_k=0$.
\end{itemize}
\end{proposition}
\begin{proof} The first part simply follows from \eqref{E:log1}, \eqref{E:log2}, Lemma \ref{L:hypebeta} and \eqref{E:invmat}. We now prove the second part. For $1\leq \ell \leq rk+1$, let $L_{\ell,m}=(L_{1,\ell,m},\dots,L_{rk+1,\ell,m})^{\tr}\in \mathbb{C}_{\infty}^{\tr}$.  Then, by \eqref{E:u1}, \eqref{E:u2} and the first part, for $1\leq \nu \leq rk+1$, we have 
\begin{equation}\label{E:normbound}
\log_q\left|\sum_{\ell=1}^{rk+1}L_{\nu,\ell,m}z_{\ell}^{q^m}\right|\leq -(rk+1)q^m.
\end{equation}
Thus, by the properties of non-archimedean norm as well as the bounds described for the entries of elements in $\mathfrak{U}$ above, we see that $\Log_{\rho^{[k]}}(\mathfrak{U})\subseteq \mathfrak{U}$. 
 
Now, as in the proof of Corollary \ref{C:ourperiod}, for $1\leq j \leq r$,  let $
\dmat_k[t]^{-n}\varpi_{k,j}= (
\mathfrak{u}_{1,n},\dots,
\mathfrak{u}_{rk+1,n})^{\tr}\in \mathbb{C}_{\infty}^{rk+1}.$ Observe, by \eqref{E:u1} and \eqref{E:u2}, that for each $n\geq 1$, $\dmat_k[t]^{-n}\varpi_{k,j}\in \mathfrak{U}$. Furthermore, we have 
  \begin{multline*}
\log_q\left|\dmat_k[t]^{-n}\varpi_{k,j}\right|\geq \log_q\left|\mathfrak{u}_{rk+1,n}\right|=\log_q\left|\frac{\tilde{\pi}^k\lambda_{j}}{\theta^{n}}\right|\\
=R_{\phi}+kR_C-n=R_{\phi}+\frac{kq}{q-1}-n>-nq^m>-(rk+1)nq^{m}
   \end{multline*}
 for all $m\geq 1$. Here, we note that  the third equality follows from \eqref{E:periodnorm}.  Thus, again by the properties of non-archimedean norm, we have
   \begin{equation}\label{E:isom}
   \left|\Log_{\rho^{[k]}}(d_k[t]^{-n}\varpi_{k,j})\right|=\left|d_k[t]^{-n}\varpi_{k,j}\right|
   \end{equation}
  as desired. Now we show that $\mathfrak{m}_k=0$. Let $\Exp_{\rho^{[k]}}=\sum_{m\geq 0}\alpha_m\tau^m\in \Mat_d(\mathbb{C}_{\infty}[[\tau]])$ with $\alpha_0=\Id_{rk+1}$. Since $\Exp_{\rho^{[k]}}\Log_{\rho^{[k]}}=\Id_{rk+1}$, for each $m\geq 0$, we have 
  \begin{equation}\label{E:expcoef}
  \alpha_{m}=-\sum_{i=1}^{m}\alpha_{m-i}L_{i}^{(m-i)}.
  \end{equation}
  Observe from the first part that 
  \[
  \log_q|L_{1,i}^{(m-i)}|\leq-q^{m}\left(R_{\phi}+\frac{k}{q-1}\right)+\mathfrak{P}_{1}
  \]
  and similarly
  \[
 \log_q|L_{ru+j+1,i}^{(m-i)}| \leq -q^m\left(q^{j-r}R_{\phi}+u+1+\frac{k}{q-1}\right)+\mathfrak{P}_{ru+j+1}, \ \ \  \ \  1\leq j \leq r, \ \ 0\leq u \leq k-1.
  \]
  Letting $\alpha_{\ell,m}$ be the $\ell$-th column of $\alpha_{m}$, by \eqref{E:expcoef}, one can then recursively obtain that 
   \[
    \log_q|\alpha_{1,m}|\leq -q^m\left(R_{\phi}+\frac{k}{q-1}\right)+\mathfrak{O}_1
    \]
    and 
    \[
    \log_q|\alpha_{ru+j+1,m}|\leq -q^m\left(q^{j-r}R_{\phi}+u+1+\frac{k}{q-1}\right)+\mathfrak{O}_{ru+j+1}, \ \ \  \ \  1\leq j \leq r, \ \ 0\leq u \leq k-1
    \]  
     for some real numbers $\mathfrak{O}_1$ and $\mathfrak{O}_{ru+j+1}$ independent of $n$.  Now, applying the same analysis used for $\Log_{\rho^{[k]}}$ above, we see that $\Exp_{\rho^{[k]}}(\mathfrak{U})\subseteq \mathfrak{U}$. Since, for each $n\geq 1$, $d_k[t]^{-n}\varpi_{k,j}\in \mathfrak{U}$, we  obtain $\Exp_{\rho^{[k]}}(d_k[t]^{-n}\varpi_{k,j})\in \mathfrak{U}$. Since $\Log_{\rho^{[k]}}$ converges on $\mathfrak{U}$, we see that $\Log_{\rho^{[k]}}$ also converges at $\Exp_{\rho^{[k]}}(d_k[t]^{-n}\varpi_{k,j})$  and hence $\mathfrak{m}_k=0$.

\end{proof}

\section{Proof of Theorem \ref{T:specval}}\label{S:specvalues}
Our goal in this section is to provide a proof for Theorem \ref{T:specval}. Our method may be seen as a generalization of the techniques used by the first author in the proof of \cite[Thm. 1.1]{Gez21}. Throughout this section, we fix a Drinfeld module $\phi$ of rank $r$ given by
\[
\phi_{t}=\theta +k_1\tau+\cdots+k_r\tau^r
\]
so that $k_1,\dots,k_{r-1}\in \mathbb{F}_q$ and $k_r\in \mathbb{F}_q^{\times}$. Moreover, we recall the abelian Anderson module $\rho^{[k]}$ defined in \S\ref{S:drinfeldtens} as well as the generators $\varpi_{k,1},\dots,\varpi_{k,r}$ for the period lattice $\Lambda_{\rho^{[k]}}$ defined in \eqref{E:periodsofext}. For each $1\leq j \leq r$, recall the element $\mathcal{L}_j\in \Mat_{1\times r}(\mathbb{T})$ defined in the statement of Theorem \ref{T:specval}. Since, in this case, by Proposition \ref{P:norms}, $\mathfrak{m}_k=0$ and by Corollary \ref{C:ourperiod}, for $i\geq 1$, each $\dmat_{k}[t]^{-i}\varpi_{k,\nu}$ lies in the domain of convergence of $\Log_{\rho^{[k]}}$, we are now able to provide an alternative proof for \cite[Thm. 1.6]{GG25} as a consequence of Proposition \ref{P:periodbasis} and the last assertion of Theorem \ref{T:mellin}.
\begin{theorem}[{cf. \cite[Thm. 1.6]{GG25}}] \label{T:tensor} Let $\zz\in \mathbb{C}_{\infty}^{rk+1}$ lie in the domain of convergence of $\Log_{\rho^{[k]}}$ and let $\mathfrak{p}_{\ell}:\mathbb{C}_{\infty}^{rk+1}\to \mathbb{C}_{\infty}$ be the projection to the $\ell$-th coordinate. For any $1\leq j \leq r$, we have
\[
\mathfrak{p}_{rk-r+j+1}(\Log_{\rho^{[k]}}(\zz))=\mathcal{M}_{\zz}\left(
\mathcal{L}_j\right) .
\]
\end{theorem}

To prove Theorem \ref{T:specval}, we start with defining \textit{the class module of $\rho^{[k]}$} which is the finite $\bA$-module given by
\[
H(\rho^{[k]}/A):=\frac{\rho^{[k]}(K_{\infty})}{\Exp_{\rho^{[k]}}(\Lie(\rho^{[k]})(K_{\infty}))+\rho^{[k]}(A)}.
\]

\begin{lemma}\label{L:classmod} We have 
\[
H(\rho^{[k]}/A)=\{0\}.
\]
\end{lemma}
\begin{proof}  Proposition \ref{P:roclog} implies that the set 
\[
\mathfrak{m}:=\{z\in \Lie(\rho^{[k]})(K_{\infty}) \ \ | \ \ |x|\leq q^{-1}\}
\]
    lies in the domain of convergence of $\Log_{\rho^{[k]}}$. Since $\Log_{\rho^{[k]}}$ is the formal inverse of $\Exp_{\rho^{[k]}}$ and we have the decomposition,
    \[
    \Lie(\rho^{[k]})(K_{\infty})=\mathfrak{m}\oplus \Lie(\rho^{[k]})(A),
    \]
    the lemma immediately follows.
\end{proof}

\begin{lemma}\label{L:linind} Let $s_i\in \Mat_{(rk+1)\times 1}(\mathbb{F}_q)$ be the $i$-th unit vector. Then for each $1\leq i \leq rk+1$, $\Log_{\rho^{[k]}}(s_i)$ converges in $\mathbb{C}_{\infty}^d$. Moreover, the set \[
\{\Log_{\rho^{[k]}}(s_1),\dots,\Log_{\rho^{[k]}}(s_{rk+1})\}
\]
is $\bA$-linearly independent. In other words, if there exist $a_1,\dots,a_{rk+1}\in \bA$ satisfying
\begin{equation}\label{E:linindid}
\dmat_k[a_1]\Log_{\rho^{[k]}}(s_1)+\cdots+ \dmat_k[a_{rk+1}]\Log_{\rho^{[k]}}(s_{rk+1})=(0,\dots,0)^{\tr},
\end{equation}
    then $a_1=\cdots=a_{rk+1}=0$.
\end{lemma}
\begin{proof} Using the functional equation \eqref{E:logfunc} and our condition on the elements, $k_1,\dots,k_r$, we see that, for all $n\geq 1$, each entry of $L_n$ lies in $\mathfrak{m}$. Thus, using the same idea to prove \cite[Prop. 4.2]{Gez21} finishes the proof of the lemma and hence we leave its details to the reader.
    
\end{proof}

Recall the unit module $U(\rho^{[k]}/A)$ of $\rho^{[k]}$ from \S\ref{S:analysisunit}.

\begin{lemma}\label{L:unit} We have 
\[
U(\rho^{[k]}/A)=\bigoplus_{i=1}^{rk+1}\bA\cdot \Log_{\rho^{[k]}}(s_i).
\]
\end{lemma}
\begin{proof} Let $\Lambda$ be the $\bA$-module generated by $ \Log_{\rho^{[k]}}(s_i)$ for $1\leq i \leq rk+1$. Then, by Lemma \ref{L:linind}, $\Lambda$ is a free $\bA$-module of rank $rk+1$. Thus, by Lemma \ref{L:classmod} and \cite[Prop. 4.10]{Gez21}, we obtain $\Lambda=U(\rho^{[k]}/A)$ as desired.
\end{proof}

Let $L(\rho^{[k]}/A)$ be the Taelman $L$-value of $\rho^{[k]}$ defined as in \S\ref{S:lfunctionsTaelman} . Since $\phi$ has everywhere good reduction, using Fang's class number formula \cite[Thm. 1.10]{Fang}, Proposition \ref{P:lvalueTaelman} and Lemma \ref{L:classmod}, we now have 
\begin{equation}\label{E:CNF}
L(\phi^{\vee},k)=L(\rho^{[k]}/A)=[\Lie(\rho^{[k]})(A):U(\rho^{[k]}/A)].
\end{equation}
In what follows, our goal is to evaluate $[\Lie(\rho^{[k]})(A):U(\rho^{[k]}/A)]$. In order to do this, we first state  a useful lemma on hyperderivatives.
\begin{lemma}\label{L:partder} Let $0\leq j \leq q-1$ and $n\in \mathbb{Z}$. We have 
\[
\partial_t^j((t-\theta^q)^n)|_{t=\theta}=\partial^j_{t}((t-t^q)^n)|_{t=\theta}.
\]
\end{lemma}
\begin{proof} We first let $n\geq 0$. 
Since 
\[
\partial_{t}^j(t-t^q)=\begin{cases}
    1 & \text{ if } j=1\\
    0 & \text{ if } 2\leq j \leq q-1,
\end{cases}
\]
by  Lemma \ref{L:hyp}(i), we obtain 
\[
\partial^j_{t}((t-t^q)^n)|_{t-\theta}=\binom{n}{j}(t-t^q)^{n-j}|_{t=\theta}=\partial^j_t((t-\theta^q)^n)|_{t=\theta}.
\]
Now, let $n=-m<0$ for some $m>0$. By Lemma \ref{L:hyp}(i), we obtain 
    \[
\partial_{t}^j((t-t^q)^n)=\partial_{t}^j\left(\frac{1}{(t-t^q)^{m}}\right)=(-1)^j\binom{m+j-1}{j}\frac{1}{(t-t^q)^{m+j}}.
    \] 
    On the other hand, by Lemma \ref{L:hyp}(ii), we have 
    \[
   \partial_{t}^j((t-\theta^q)^n)=\partial^j_{t}\left(\frac{1}{(t-\theta^q)^m}\right)=(-1)^j\binom{m+j-1}{j}\frac{1}{(t-\theta^q)^{m+j}}. 
    \]
Above calculations, after evaluating at $t=\theta$, finish the proof of the lemma.
    
\end{proof}

Next, inspired by  \cite[Prop. 4.2]{ANDTR20}, we construct a certain $\mathbb{F}_q((1/t))$-vector subspace of $\Lie(\rho^{[k]})(K_{\infty})$ which will be later useful to prove our result.

\begin{proposition}\label{P:intersection} Let $1\leq k\leq q-1$. For $1\leq i \leq r$, set 
\[
\mathfrak{h}_i:=\widehat{\delta_0\circ\iota}\left(\mathfrak{e}_{i}(\mathfrak{J}^{-1})^{(1)}\right)\in \Mat_{(rk+1)\times 1}(K_{\infty}).
\]
    Then, the set $\{\mathfrak{h}_1,\dots,\mathfrak{h}_r\}$ is $\mathbb{F}_q((1/t))$-linearly independent. Moreover, 
    if we let
\[
\mathfrak{W}:=\bigoplus_{i=1}^{r}\mathbb{F}_q((1/t))\cdot \mathfrak{h}_i,
\]
then
    \[
\mathfrak{W}=\bigoplus_{i=rk-r+2}^{rk+1}\mathbb{F}_q((1/t))\cdot s_{i}.
    \]
\end{proposition}
\begin{proof} Note that, we have
\[
(t-\theta^q)^{k}(\mathfrak{J}^{-1})^{(1)}=\begin{pmatrix}
 & \mathfrak{w}_{12} &\cdots & \mathfrak{w}_{1r}\\
 & & \ddots & \vdots\\
 & & & \mathfrak{w}_{(r-1)r}\\
 \mathfrak{w}_{r1}(t-\theta^{q})^{-1}&\mathfrak{w}_{r2} & \cdots & \mathfrak{w}_{rr}
\end{pmatrix}
\]
where the entries in the first $r-1$-row of the right hand side above are in $\mathbb{F}_q$, $\mathfrak{w}_{12},\dots,\mathfrak{w}_{(r-2)(r-1)}$, $\mathfrak{w}_{(r-1)r},\mathfrak{w}_{r1}\in \mathbb{F}_q^{\times}$ and $\mathfrak{w}_{r2},\dots,\mathfrak{w}_{rr}\in K(t)$. We have
\begin{multline*}
\mathfrak{h}_{r-1}=\widehat{\delta_0\circ\iota}\left(\left[0, \dots,0, \frac{\mathfrak{w}_{(r-1)r}}{(t-\theta^q)^{k}}\right]\right)
= [\partial_t^k\left(\frac{\mathfrak{w}_{(r-1)r}}{(t-\theta^q)^{k}}\right),0, \dots,0,\\
\partial_t^{k-1}\left(\frac{\mathfrak{w}_{(r-1)r}}{(t-\theta^q)^{k}}\right),0,\dots, ,0,\partial_t^1\left(\frac{\mathfrak{w}_{1r}}{(t-\theta^q)^{k}}\right), 0,\dots,0,\frac{\mathfrak{w}_{(r-1)r}}{(t-\theta^q)^{k}} ]^{\tr} |_{t=\theta}.
\end{multline*}
To be more precise, we remark that, above, for $0\leq \nu \leq k$, $\partial_t^{\nu}\left(\frac{\mathfrak{w}_{(r-1)r}}{(t-\theta^q)^{k}}\right)$ appears in the $(r\nu+1)$-th entry of $\mathfrak{h}_{r-1}$. On the other hand, using the Leibniz rule for hyperderivatives and Lemma \ref{L:partder}, for any $0\leq m \leq k\leq q-1$ and $\mu\geq 1$, we obtain
\begin{multline}\label{E:h1}
		\sum_{\substack{\nu.j\geq 0\\        \nu+j=m}}\partial_{t}^j((\theta-\theta^q)^{\mu})|_{t=\theta}\partial_{t}^{\nu}\Big(\frac{1}{(t-\theta^q)^{\mu}}\Big)_{|t=\theta}=\sum_{\substack{\nu.j\geq 0\\\nu+j=m}}\partial_{t}^j((t-\theta^q)^{\mu})_{|t=\theta}\partial_{t}^{\nu}\Big(\frac{1}{(t-\theta^q)^{\mu}}\Big)_{|t=\theta}=\partial_{t}^m(1).
\end{multline}
Observe that the left hand side of \eqref{E:h1}, after multiplying with $\mathfrak{w}_{(r-1)r}$ and letting $\mu=k$, is indeed the $ (r(k-m)+1)$-th entry of $\dmat_k[(\theta-\theta^q)^k]\mathfrak{h}_{r-1}$. Therefore, we have 
\begin{equation}\label{E:linearity0}
\mathfrak{w}_{(r-1)r}^{-1}(\theta-\theta^q)^k\cdot \mathfrak{h}_{r-1}=\mathfrak{w}_{(r-1)r}^{-1}\dmat_k[(\theta-\theta^q)^k]\mathfrak{h}_{r-1}=[0,\dots,0,1]=s_{rk+1}.
\end{equation}

We now apply the same idea for $\mathfrak{h}_{r-2}$ to obtain $s_{rk}$ in terms of a linear combination of $\mathfrak{h}_{r-1}$ and $\mathfrak{h}_{r-2}$. Observe that we have 
\[
\mathfrak{h}_{r-2}=\widehat{\delta_0\circ\iota}\left(\left[0, \dots,0, \frac{\mathfrak{w}_{(r-2)(r-1)}}{(t-\theta^q)^{k+1}},\frac{\mathfrak{w}_{(r-2)r}}{(t-\theta^q)^{k+1}}\right]\right).
\]
Therefore, we see that for $0\leq \nu \leq k$, $\partial_t^{\nu}\left(\frac{\mathfrak{w}_{(r-2)r}}{(t-\theta^q)^{k}}\right)$ appears in the $(r\nu+1)$-th entry of $\mathfrak{h}_{r-2}$ and $\partial_t^{\nu}\left(\frac{\mathfrak{w}_{(r-2)(r-1)}}{(t-\theta^q)^{k}}\right)$ appears in the $(r\nu+2)$-th entry of $\mathfrak{h}_{r-2}$. Using \eqref{E:h1}, we now obtain that 
\[
\mathfrak{w}_{(r-2)(r-1)}^{-1}(\theta-\theta^q)^{k}\cdot \mathfrak{h}_{r-2}-\mathfrak{w}_{(r-2)r}\mathfrak{w}_{(r-1)r}^{-1}(\theta-\theta^q)^k\cdot \mathfrak{h}_{r-1}=[0,\dots,1,0]=s_{rk}.
\]
Proceeding in this way, we see that for $1\leq i \leq r-2$, there exist $c_{r-i},\dots,c_{r-1}\in \mathbb{F}_q$ with $c_{r-i}\in \mathbb{F}_q^{\times}$ such that 
\[
c_i(\theta-\theta^q)^k\cdot \mathfrak{h}_{r-i}+\cdots+c_{r-1}(\theta-\theta^q)^k\cdot \mathfrak{h}_{r-1}=s_{rk+1-i}.
\]
Similarly, there exist $c_1',\dots,c_{r}'\in \mathbb{F}_q$ with $c_r'\in \mathbb{F}_q^{\times}$ such that 
\begin{equation}\label{E:linearity}
c_1'(\theta-\theta^q)^{k+1}\cdot\mathfrak{h}_1+\cdots+c_{r}'(\theta-\theta^q)^{k+1}\cdot\mathfrak{h}_{r} =s_{rk-r+2}.
\end{equation}
Thus, \eqref{E:linearity0} and \eqref{E:linearity} imply 
\[
\bigoplus_{i=rk-r+2}^{rk+1} \mathbb{F}_q((1/t))\cdot s_i \subseteq \sum_{i=1}^r \mathbb{F}_q((1/t))\cdot \mathfrak{h}_i.
\]
Since the $\mathbb{F}_q((1/t))$-vector space in the left hand side above has dimension $r$, the set  $\{\mathfrak{h}_1,\dots,\mathfrak{h}_r\}$ must be $\mathbb{F}_q((1/t))$-linearly independent. Thus, we have 
\[
\mathfrak{W}=\bigoplus_{i=1}^{r}\mathbb{F}_q((1/t))\cdot \mathfrak{h}_i=\bigoplus_{i=rk-r+2}^{rk+1}\mathbb{F}_q((1/t))\cdot s_{i}
\]
as desired.
\end{proof}
\begin{proposition}\label{P:fin} Let $1\leq k \leq q-1$. 
     For any $rk-r+2\leq i \leq rk+1$, we have 
$
\Log_{\rho^{[k]}}(s_i)\in \mathfrak{W}.
$   Moreover 
\[
U(\rho^{[k]}/A)\cap \mathfrak{W}=\bigoplus_{i=rk-r+2}^{rk+1}\bA\cdot \Log_{\rho^{[k]}}(s_i) \ \ \text{ and } \ \ 
\Lie(\rho^{[k]})(A)\cap \mathfrak{W}=\bigoplus_{i=rk-r+2}^{rk+1}\bA\cdot s_i.
\]
Furthermore, $U(\rho^{[k]}/A)\cap \mathfrak{W}$ and $\Lie(\rho^{[k]})(A)\cap \mathfrak{W}$ are $\bA$-lattices in $\mathfrak{W}$.
\end{proposition}
\begin{proof} Our first goal is to show that 
\[
\Log_{\rho^{[k]}}(s_i)\in \mathfrak{W}=\bigoplus_{i=rk-r+2}^{rk+1} \mathbb{F}_q((1/t))\cdot s_i. 
\]
Observe that, for $n\geq 0$, we may write 
\begin{multline*}
\mathfrak{J}^{[n]}=\frac{1}{\mathbb{S}_n^{k+1}}\times \\
\begin{pmatrix}
    \displaystyle\sum_{0\leq i_1,\dots,i_n\leq 1} \mathfrak{w}_{1,1,i_1,\dots,i_n}(t-\theta^q)^{i_1}\cdots (t-\theta^{q^n})^{i_n}& \cdots & \displaystyle\sum_{0\leq i_1,\dots,i_n\leq 1}\mathfrak{w}_{1,r,i_1,\dots,i_n}(t-\theta^q)^{i_1}\cdots (t-\theta^{q^n})^{i_n}\\
    \vdots & & \vdots \\
   \displaystyle\sum_{0\leq i_1,\dots,i_n\leq 1} \mathfrak{w}_{r,1,i_1,\dots,i_n}(t-\theta^q)^{i_1}\cdots (t-\theta^{q^n})^{i_n}&\cdots  & \displaystyle\sum_{0\leq i_1,\dots,i_n\leq 1} \mathfrak{w}_{r,r,i_1,\dots,i_n}(t-\theta^q)^{i_1}\cdots (t-\theta^{q^n})^{i_n}\end{pmatrix}
\end{multline*}
    for some $\mathfrak{w}_{\nu,\mu,i_1,\dots,i_n}\in \mathbb{F}_q$ where $1\leq \nu,\mu\leq r$. Recall that $\mathfrak{J}_i^{[n]}$ is the $i$-th row of $\mathfrak{J}^{[n]}$ and let $L_{\ell,n}$ be the $\ell$-th column of the $n$-th logarithm coefficient of $\Log_{\rho^{[k]}}$. Then by \eqref{E:log1} and \eqref{E:log2}, for $0\leq u \leq k-1$ and $1\leq j \leq r$, we have 
    \[
    L_{1,n}=
    \widehat{\delta_0\circ \iota}\left(\frac{(t-\theta^{q^n})^{k}}{\mathbb{S}_n^{k+1}}\mathfrak{J}^{[n]}_{r}\right)   
   \]
    and 
    \[
    L_{ru+j+1}=  \widehat{\delta_0\circ \iota}\left(\frac{(t-\theta^{q^n})^{k-u-1}}{\mathbb{S}_n^{k+1}}\mathfrak{J}^{[n]}_{j}\right).    \]
    By Lemma \ref{L:partder} and the Leibniz rule, for $0\leq u \leq q-1$ and $1\leq \nu,\mu\leq r$, we see that 
    \begin{multline}\label{E:logcoef}
    \partial_t^{u}\left(\frac{(t-\theta^{q^n})^{k-u}}{\mathbb{S}_n^{k+1}}\displaystyle\sum_{0\leq i_1,\dots,i_n\leq 1}\mathfrak{w}_{\nu,\mu,i_1,\dots,i_n}(t-\theta^q)^{i_1}\cdots (t-\theta^{q^n})^{i_n} \right)|_{t=\theta}\\
    =\partial_{t}^{u}\left(\frac{(t-t^{q^n})^{k-u}}{((t-t^q)\cdots (t-t^{q^n}))^{k+1}}\displaystyle\sum_{0\leq i_1,\dots,i_n\leq 1}\mathfrak{w}_{\nu,\mu,i_1,\dots,i_n}(t-t^q)^{i_1}\cdots (t-t^{q^n})^{i_n} \right)|_{t=\theta}.\end{multline}
    Now, for $n\geq 0$ and $1\leq \mu \leq r$, let
\[
\mathfrak{L}_{1,\mu,n}:=\frac{(t-t^{q^n})^{k}}{((t-t^q)\cdots (t-t^{q^n}))^{k+1}}\displaystyle\sum_{0\leq i_1,\dots,i_n\leq 1}\mathfrak{w}_{r,\mu,i_1,\dots,i_n}(t-t^q)^{i_1}\cdots (t-t^{q^n})^{i_n}\in \mathbb{F}_q(t)
\]
and 
\[
\mathfrak{L}_{ru+j+1,\mu,n}:=\frac{(t-t^{q^n})^{k-u-1}}{((t-t^q)\cdots (t-t^{q^n}))^{k+1}}\displaystyle\sum_{0\leq i_1,\dots,i_n\leq 1}\mathfrak{w}_{j,\mu,i_1,\dots,i_n}(t-t^q)^{i_1}\cdots (t-t^{q^n})^{i_n}\in \mathbb{F}_q(t).
\]
Then, \eqref{E:logcoef} implies that 
\begin{equation}\label{E:logcof1}
\dmat_k[\mathfrak{L}_{1,1,n}] s_{rk-r+2}+\cdots+\dmat_k[\mathfrak{L}_{1,r,n}] s_{rk+1}=L_{1,n}\in \mathfrak{W}
\end{equation}
and for $1\leq j \leq r$, $0\leq u \leq k-1$, we have
\begin{equation}\label{E:logcof2}
\dmat_k[\mathfrak{L}_{ru+j+1,1,n}] s_{rk-r+2}+\cdots+\dmat_k[\mathfrak{L}_{ru+j+1,r,n}] s_{rk+1}=L_{ru+j+1,n}\in \mathfrak{W}.
\end{equation}
Thus, since for $n\geq 0$ and $ rk-r+2\leq i \leq rk+1 $, each $L_{i,n}s_{i}\in \mathfrak{W}$ and $\mathfrak{W}$ is a complete normed $\mathbb{F}_q((1/t))$-vector space, we have 
\[
\Log_{\rho^{[k]}}(s_{i})=\sum_{n=0}^{\infty}L_{i,n}s_{i}\in \mathfrak{W}.
\]
We now prove the second assertion. By Lemma \ref{L:unit} and the first assertion, note that 
\[
U(\rho^{[k]}/A)\cap \mathfrak{W}\supseteq \sum_{i=rk-r+2}^{rk+1}\bA\cdot \Log_{\rho^{[k]}}(s_i)= \bigoplus_{i=rk-r+2}^{rk+1}\bA\cdot \Log_{\rho^{[k]}}(s_i)
\]
where the last equality follows from Lemma \ref{L:linind}. On the other hand, since $\dim_{\mathbb{F}_q((1/t))}(U(\rho^{[k]}/A)\cap \mathfrak{W})\leq  \dim_{\mathbb{F}_q((1/t))}(\mathfrak{W})=r$, we obtain 
\[
U(\rho^{[k]}/A)\cap \mathfrak{W}\subseteq \bigoplus_{i=rk-r+2}^{rk+1}\bA\cdot \Log_{\rho^{[k]}}(s_i),
\]
finishing the proof of the second assertion. The third assertion can be obtained similarly and we leave its details to the reader. Finally, the last assertion follows from \cite[Thm. 4.4]{ANDTR20}.
\end{proof}

\begin{proof}[{Proof of Theorem \ref{T:specval}}] 
By \cite[\S2.1]{AR16}, note that  $ \frac{\Lie(\rho^{[k]})(A)}{\Lie(\rho^{[k]})(A)\cap \mathfrak{W}} $ and $ \frac{U(\rho^{[k]}/A)}{U(\rho^{[k]}/A)\cap \mathfrak{W}} $ are $\bA$-lattices in $\frac{\Lie(\rho^{[k]})(K_{\infty})}{\mathfrak{W}}$. Moreover, we have 
\begin{multline*}
L(\phi^{\vee},k)=[\Lie(\rho^{[k]})(A):U(\rho^{[k]}/A)]_{|t=\theta}\\
=[\Lie(\rho^{[k]})(A)\cap \mathfrak{W}:U(\rho^{[k]}/A)\cap \mathfrak{W}]_{|t=\theta}\left[\frac{\Lie(\rho^{[k]})(A)}{\Lie(\rho^{[k]})(A)\cap \mathfrak{W}}:\frac{U(\rho^{[k]}/A)}{U(\rho^{[k]}/A)\cap \mathfrak{W}}\right]|_{t=\theta}
\end{multline*}
where the last equality follows from \cite[\S2.1]{AR16}. 
We first note that, since $\mathfrak{W}$ is a complete normed $\mathbb{F}_q((1/t))$-vector space, by \eqref{E:logcof1} and \eqref{E:logcof2}, we have 
\[
\sum_{n=1}^{\infty}L_{1,n}=\Log_{\rho^{[k]}}(s_{1})-s_{1}\in \mathfrak{W}
\]
and 
\[
\sum_{n=1}^{\infty}L_{ru+j+1,n}=\Log_{\rho^{[k]}}(s_{ru+j+1})-s_{ru+j+1}\in \mathfrak{W}, \ \ \ \\ \ \ 1\leq j \leq r, \ \ 0\leq u \leq k-1.
\]
In particular, for $1\leq i \leq rk-1$, we have
$
\Log_{\rho^{[k]}}(s_{i})\equiv s_{i} \mod \mathfrak{W}.
$
Therefore the $\bA$-module $ \frac{U(\rho^{[k]}/A)}{U(\rho^{[k]}/A)\cap \mathfrak{W}}\cong \bigoplus_{i=1}^{rk-r+1}\bA\cdot \Log_{\rho^{[k]}}(s_i) $ is generated by $s_{1},\dots,s_{rk-r+1}$ in $\frac{\Lie(\rho^{[k]})(K_{\infty})}{ \mathfrak{W}}$. Similarly, the $\bA$-module $\frac{\Lie(\rho^{[k]})(A)}{\mathfrak{W}} $ is generated by $s_{1},\dots,s_{rk-r+1}$ in $\frac{\Lie(\rho^{[k]})(K_{\infty})}{\mathfrak{W}}$. 
Thus, we obtain
\[
\left[\frac{\Lie(\rho^{[k]})(A)}{\Lie(\rho^{[k]})(A)\cap \mathfrak{W}}:\frac{U(\rho^{[k]}/A)}{U(\rho^{[k]}/A)\cap \mathfrak{W}}\right]|_{t=\theta}=\det(\Id_{r(k-1)+1})=1.
\]
Therefore, we have
\[
L(\phi^{\vee},k)=[\Lie(\rho^{[k]})(A)\cap \mathfrak{W}:U(\rho^{[k]}/A)\cap \mathfrak{W}]|_{t=\theta}.
\]
Finally, using the same idea in \eqref{E:tractablereg}, we see that $[\Lie(\rho^{[k]})(A)\cap \mathfrak{W}:U(\rho^{[k]}/A)\cap \mathfrak{W}]|_{t=\theta} $ is the determinant of the $r\times r$-matrix whose $i$-th row is the last $r$ entries of $\Log_{\rho^{[k]}}(s_{r(k-1)+i+1})$. Since, this determinant is monic, it is equal to $L(\phi^{\vee},k)$ and thus, by Theorem \ref{T:tensor}, we finish the proof of the theorem.
\end{proof}

\end{document}